\documentclass{amsart}
\usepackage{amssymb, amscd, latexsym, graphicx, psfrag, multirow}

\newtheorem{dummy}{dummy}[section]

\newtheorem{theorem}[dummy]{Theorem}
\newtheorem{conjecture}[dummy]{Conjecture}
\newtheorem{corollary}[dummy]{Corollary}
\newtheorem{proposition}[dummy]{Proposition}

\theoremstyle{definition}

\theoremstyle{remark}
\newtheorem{remark}[dummy]{Remark}

\newcommand{\Si}{\Sigma}
\newcommand{\Ga}{\Gamma}

\newcommand{\si}{\sigma}

\newcommand{\bm}{\mathbf{m}}
\newcommand{\bp}{\mathbf{p}}
\newcommand{\bq}{\mathbf{q}}

\newcommand{\bae}{\bar{e}}

\newcommand{\cA}{\mathcal{A}}
\newcommand{\cC}{\mathcal{C}}
\newcommand{\cD}{\mathcal{D}}
\newcommand{\cM}{\mathcal{M}}
\newcommand{\cO}{\mathcal{O}}
\newcommand{\cP}{\mathcal{P}}

\newcommand{\bC}{ \mathbb{C} }
\newcommand{\bF}{ \mathbb{F} }

\newcommand{\bL}{ {\mathbb{L}} }
\newcommand{\bP}{ \mathbb{P} }
\newcommand{\bQ}{ \mathbb{Q} }
\newcommand{\bR}{ \mathbb{R} }
\newcommand{\bT}{ \mathbb{T} }
\newcommand{\bZ}{ \mathbb{Z} }

\newcommand{\Aut}{\mathrm{Aut}}
\newcommand{\Def}{\mathrm{Def}}
\newcommand{\Hom}{\mathrm{Hom}}
\newcommand{\Ker}{\mathrm{Ker}}
\newcommand{\Obs}{\mathrm{Obs}}
\newcommand{\Spec}{\mathrm{Spec}}

\newcommand{\eff}{\mathrm{eff}}

\newcommand{\rI}{\mathrm{I}}
\newcommand{\rII}{\mathrm{II}}

\newcommand{\vir}{ {\mathrm{vir}}  }
\newcommand{\ev}{\mathrm{ev}}

\newcommand{\fg}{\mathfrak{g}}
\newcommand{\ft}{\mathfrak{t}}

\newcommand{\wf}{\mathsf{f}}
\newcommand{\wu}{\mathsf{u}}
\newcommand{\wv}{\mathsf{v}}

\newcommand{\ww}{\mathsf{w}}
\newcommand{\wx}{{x}}
\newcommand{\wy}{{y}}

\newcommand{\tC}{\widetilde{C}}

\newcommand{\tI}{\widetilde{I}}
\newcommand{\tM}{\widetilde{M}}
\newcommand{\tN}{\widetilde{N}}
\newcommand{\tX}{\widetilde{X}}
\newcommand{\tbeta}{ {{\beta'}} }
\newcommand{\tbT}{\widetilde{\bT}}

\newcommand{\tbL}{ {\widetilde{\bL}}}

\newcommand{\tl}{\widetilde{l}}
\newcommand{\tq}{{q}}
\newcommand{\tv}{\widetilde{v}}
\newcommand{\tx}{\widetilde{x}}
\newcommand{\ty}{\widetilde{y}}

\newcommand{\tmu}{\widetilde{\mu}}

\newcommand{\hB}{\hat{B}}
\newcommand{\hD}{\hat{D}}
\newcommand{\hL}{\hat{L}}
\newcommand{\hT}{\hat{T}}
\newcommand{\hX}{\hat{X}}
\newcommand{\hY}{\hat{Y}}

\newcommand{\hu}{\hat{u}}

\newcommand{\vd}{\vec{d}}
\newcommand{\vw}{\vec{w}}
\newcommand{\vx}{{\vec{x}}}
\newcommand{\vmu}{\vec{\mu}}
\newcommand{\vnu}{\vec{\nu}}

\renewcommand{\hbar}{{z}}
\newcommand{\Mbar}{\overline{\cM}}
\newcommand{\bSi}{\partial\Si}
\newcommand{\lra}{\longrightarrow}
\newcommand{\tri}{\triangle}
\newcommand{\bu}{\bullet}

\begin{document}

\title{Open Gromov-Witten Invariants of Toric Calabi-Yau 3-Folds}

\author{Bohan Fang}
\address{Bohan Fang, Department of Mathematics, Columbia University,
2990 Broadway, New York, NY 10027}
\email{b-fang@math.columbia.edu}

\author{Chiu-Chu Melissa Liu}
\address{Chiu-Chu Melissa Liu, Department of Mathematics, Columbia University,
2990 Broadway, New York, NY 10027} \email{ccliu@math.columbia.edu}

\begin{abstract}
We present a proof of the mirror conjecture of Aganagic-Vafa \cite{AV} and 
Aganagic-Klemm-Vafa \cite{AKV} on disk enumeration in toric Calabi-Yau 3-folds for
all smooth semi-projective toric Calabi-Yau 3-folds. 
We consider both inner and outer branes, at arbitrary framing.  
In particular, we recover previous results on the conjecture for 
(i) an inner brane at zero framing in $K_{\bP^2}$ (Graber-Zaslow \cite{GZ}), 
(ii) an outer brane at arbitrary framing in the resolved conifold
$\cO_{\bP^1}(-1)\oplus \cO_{\bP^1}(-1)$ (Zhou \cite{Zh3}), and (iii) an outer brane at
zero framing in $K_{\bP^2}$ (Brini \cite{Br}).
\end{abstract}

\maketitle

\tableofcontents

\section{Introduction}

\subsection{Open Gromov-Witten invariants}
\label{sec:openGw}

Let $X$ be a K\"{a}hler manifold
(or more generally, an almost K\"{a}hler manifold) and
let $L$ be a Lagrangian submanifold of $X$. Naively, 
open Gromov-Witten (GW) invariants of the pair $(X,L)$ count
holomorphic maps (or more generally, $J$-holomorphic maps)
from bordered Riemann surfaces to $X$ with boundaries
mapped into $L$. Open GW theory can be viewed as a 
mathematical theory of A-model topological open strings. String
dualities have produced many striking predictions on open GW invariants, 
especially when $X$ is a Calabi-Yau 3-fold and $L$ is a special Lagrangian
submanifold (or more generally, when the Maslov class of $L$ is zero).

In this paper, we consider toric Calabi-Yau 3-folds
and a particular class of Lagrangian submanifolds
(called Aganagic-Vafa A-branes in this paper)
introduced by Aganagic-Vafa in \cite{AV}. 
Let $X$ be a toric Calabi-Yau 3-fold and $L$ be an Aganagic-Vafa
A-brane in $X$ (in particular, $L$ is diffeomorphic to $S^1\times\bR^2$). 
The open GW invariants of the pair $(X,L)$ depend on the following discrete data: 
\begin{itemize}
\item topological type $(g,h)$ of the domain bordered Riemann surface $\Si$, 
where $g$ is the number of handles and $h$ is the number of holes
when $\Si$ is a smooth bordered Riemann surface;
\item topological type of the map $u:(\Si,\bSi =\coprod_{i=1}^h R_i )\to (X,L)$,
which is given by the degree $\beta' =u_*[\Si]\in H_2(X,L;\bZ)$
and winding numbers $\mu_i = u_*[R_i]\in H_1(L;\bZ) =\bZ$;
\item framing $f\in \bZ$ of the Aganagic-Vafa A-brane $L$. 
\end{itemize}
More generally, $L$ can be a disjoint union
of framed A-branes $(L_1,f_1), \ldots, (L_s,f_s)$,
and the open GW invariants of $(X,L)$ are
\begin{equation}\label{eqn:openGW}
N^{X,L}_{g,\beta', \mu^1,\ldots,\mu^s}(f_1,\ldots,f_s) \in \bQ
\end{equation}
where $\mu^i=(\mu^i_1,\ldots, \mu^i_{h_i})$ are the winding
numbers associated to $L_i$, and $h=\sum_{i=1}^s h_i$.

When $L_1,\ldots,L_s$ are outer branes, 
$\mu^i_j$ are positive integers, and
$\mu^i$ can be viewed as a (possibly empty) partition.
In this case, the authors of \cite{LLLZ} provided a mathematical 
definition of the open GW invariants in \eqref{eqn:openGW} 
as formal relative GW invariants  of a relative formal toric 
Calabi-Yau (FTCY) 3-fold $(\hY,\hD)$\footnote{In this paper, we use slightly
different convention. See \eqref{eqn:XLYD}.}
\begin{equation}\label{eqn:relativeGW}
N^{X,L}_{g,\beta',\mu^1,\ldots,\mu^s} (f_1,\ldots,f_s) :=
\sum_{\pi(\vd')=\beta'} F^{\hY,\hD}_{g,\vd',\mu^1,\ldots,\mu^s}. 
\end{equation}
On the right hand side of \eqref{eqn:relativeGW}, $\pi$ is a surjective
group homomorphism $H_2(\hY;\bZ)\to H_2(X,L;\bZ)$, and the relative FTCY 3-fold
$(\hY,\hD)$ is determined by the toric Calabi-Yau 3-fold $X$
and the {\em framed} Aganagic-Vafa A-branes $(L_1,f_1),\ldots,
(L_s,f_s)$. More details will be given
in Section \ref{sec:openGWouter}.

\subsection{Large $N$ duality: conjectures and results} 
\label{sec:largeN}
The large $N$ duality relates A-model topological
string theory on Calabi-Yau 3-folds to Chern-Simons gauge
theory on 3-manifolds. Motivated by the large $N$ duality, 
Aganagic-Klemm-Mari\~{n}o-Vafa proposed the topological
vertex \cite{AKMV}, an algorithm of computing 
generating functions 
\begin{equation}\label{eqn:Fbeta}
\begin{aligned}
& F_\beta^X(\lambda) = \sum_g N^X_{g,\beta} \lambda^{2g-2}, \\
& F_{\beta',\mu^1,\ldots,\mu^s}^{X,L}(\lambda;f_1,\ldots,f_s) = 
\sum_g N^{X,L}_{g,\beta',\mu^1,\ldots,\mu^s}(f_1,\ldots,f_s) \lambda^{2g-2+h},
\end{aligned}
\end{equation} 
where $\beta\in H_2(X;\bZ)$ and $N^X_{g,\beta}$ is the genus $g$, degree $\beta$
(closed) GW invariants of $X$.  Note that the generating functions
in \eqref{eqn:Fbeta} are obtained by fixing a topological type of the map and
summing over all possible topological types of the domain.

Let $Z_\beta^X$ and $Z_{\beta',\mu^1,\ldots, \mu^s}^{X,L}$ be the disconnected versions
of $F_\beta^X$ and $F_{\beta',\mu^1,\ldots,\mu^s}^{X,L}$ respectively. 
The algorithm of AKMV consists of two steps: 
\begin{enumerate}
\item[O1.] Explicit gluing formula which expresses $Z_\beta^X$, $Z_{\beta',\mu^1,\ldots, \mu^s}^{X,L}$ in terms
of the topological vertex $C_{\mu^1,\mu^2,\mu^3}$, a generating function of open 
GW invariants of $\bC^3$ relative to three Aganagic-Vafa A-branes.
\item[O2.] Evaluation of the topological vertex $C_{\mu^1,\mu^2,\mu^3}$ by 
relating it to  the colored  HOMFLY polynomial of a 3-component link in $S^3$.
\end{enumerate} 

In \cite{LLLZ}, Li-Liu-Liu-Zhou developed a mathematical theory
of the topological vertex using algebraic relative GW theory developed by 
J. Li \cite{Li1,Li2}\footnote{The symplectic relative GW theory was 
developed independently by Li-Ruan \cite{LR} and Ionel-Parker \cite{IP1, IP2}.}:
\begin{enumerate}
\item[R0.] Open GW invariants of $(X,L)$ are defined  
as formal relative GW invariants of a relative FTCY 3-fold
$(\hY,\hD)$. 
\item[R1.] The degeneration formula  satisfied by these
formal relative GW invariants agrees with the gluing formula in \cite{AKMV},
with $C_{\mu^1,\mu^2,\mu^3}$ replaced by a generating function 
$\tC_{\mu^1,\mu^2,\mu^3}$ of invariants of an in-decomposable relative FTCY 3-fold.
\item[R2.] LLLZ evaluated $\tC_{\mu^1,\mu^2,\mu^3}$  using virtual localization and a formula 
of Hodge integrals provided in  \cite{LLZ2}, and showed that it agrees with the formula of 
$C_{\mu^1,\mu^2,\mu^3}$  in \cite{AKMV}  when one of the partitions is empty. 
\end{enumerate}
This proves the topological vertex algorithm up  to 2-leg vertex. The validity of the full 3-leg case is a 
consequence of the proof of Gromov-Witten/Donaldson-Thomas correspondence of toric Calabi-Yau 3-folds 
by Maulik-Okounkov-Oblomkov-Pandharipande \cite{MOOP}.

\subsection{Mirror Symmetry: conjectures and results} \label{sec:MS}
The mirror symmetry relates A-model topological
string theory on a Calabi-Yau 3-fold $X$  to the B-model topological string
theory on the mirror Calabi-Yau 3-fold $X^\vee$. 

Let $X$ be a smooth toric Calabi-Yau 3-fold. We assume that $X$ is {\em semi-projective}, i.e.,
$X$ has a torus fixed point and $X$ is projective over its affinization $\mathrm{Spec}\big(H^0(X,\cO_X)\big)$.
By the results in \cite[Section 2]{HS}, a smooth toric variety is semi-projective if and
only if it is equal to the GIT quotient of an affine space $\bC^r$ by the action of
a subtorus of $(\bC^*)^r$.  Let $L$ be an Aganagic-Vafa A-brane in $X$. 
Aganagic-Vafa related a generating function of disk invariants of the pair $(X,L)$ 
to the Abel-Jacobi map of the mirror curve
of $X$ \cite{AV}.  It was clarified in \cite{AKV} that the framing
of $L$ corresponds to choice of flat coordinates in 
the B-model. The integrals in \cite{AV, AKV} are solutions to extended 
Picard-Fuchs equations \cite{LM, LMW, M, M2}. 
To our knowledge, the above mirror conjectures on disk invariants has been verified 
in the following cases:
(i) $X$ is the total space of $K_{\bP^2}$ and $L$ is
an inner brane at zero framing (T. Graber and E. Zaslow \cite{GZ});
(ii)$X$ is the resolved conifold (i.e. the total space of $\cO_{\bP^1}(-1)\oplus \cO_{\bP^1}(-1)$) 
and $L$ is an outer brane at arbitrary framing  (J. Zhou \cite{Zh3});
(iii) $X$ is the total space of $K_{\bP^2}$ and $L$ is 
an outer brane at zero framing (A. Brini \cite[Section 5.3]{Br}).

Based on the work of Eynard-Orantin \cite{EO} and Mari\~{n}o \cite{Ma}, 
Bouchard-Klemm-Mari\~{n}o-Pasquetti \cite{BKMP} proposed
a new formalism of the B-model topological strings on the mirrors of
toric Calabi-Yau 3-folds, and provided an recursive algorithm,
now known as the remodeling conjecture, 
which determines the generating functions
\begin{equation}\label{eqn:Fg}
\begin{aligned}
& F^X_g(Q) = \sum_{\beta\in H_2(X;\bZ)} N^X_{g,\beta} Q^\beta \\
& F^{X,L}_{g,h}(Q, x_1,\ldots,x_h;f) = 
\sum_{\tiny \begin{array}{c}\beta'\in H_2(X,L;\bZ)\\
\mu=(\mu_1,\ldots,\mu_h)\end{array}} N_{\beta',\mu_1,\ldots,\mu_h}(f) Q^{\beta'} x_1^{\mu_1}\cdots x_h^{\mu_h}
\end{aligned}
\end{equation}
from the disk invariants $F^{X,L}_{0,1}$ and annulus invariants $F^{X,L}_{0,2}$.
Note that the generating functions are obtained by fixing the topological type of the domain
and summing over all possible topological types of the map.  The remodeling conjecture 
has been proved for the framed 1-leg topological vertex by L. Chen \cite{Ch} and 
by J. Zhou \cite{Zh1}. J. Zhou later proved the conjecture
for the framed 3-leg topological vertex \cite{Zh2}. Recently,
Bouchard-Catuneanu-Marchal-Su\l kowski proved the
remodeling conjecture for $F_g^{\bC^3}$ \cite{BCOS}. 

\subsection{Summary of results} 
The goals of the present paper are twofold. 
\begin{enumerate}
\item[(i)] We use localization to derive a formula of formal 
relative GW invariants of the relative FTCY 3-fold $(\hY,\hD)$ in terms of descendant GW invariants of
the FTCY 3-fold $\hX:=\hY-\hD$ ({\bf Proposition \ref{pro:YDX}}). As a consequence, we obtain
a formula of open GW invariants of $(X,L)$ in terms of descendant GW invariants of $X$
 when $L$ is a union of framed outer branes ({\bf Proposition \ref{pro:XLX}}).
The formula allows us to extend the definition to a single framed inner brane
({\bf Proposition \ref{prop:open-closed}}).

\item[(ii)] We use the approach in Graber-Zaslow \cite{GZ} to
prove the mirror conjecture on disk invariants in   \cite{AV, AKV, LM, LMW, M, M2}
for any smooth semi-projective  toric Calabi-Yau 3-folds. 
We consider both outer and inner branes, at arbitrary framing ({\bf Conjecture \ref{AVconj}}).
We list explicit mirror formulae for the resolved conifold, local toric Fano surfaces, and
toric crepant resolutions of $\big(\cO_{\bP^1}(-1)\oplus \cO_{\bP^1}(-1)\big)/\bZ_m$.
\end{enumerate}


\subsection{Outline of the paper}
In Section \ref{sec:geometry}, we review the
geometry and topology of toric Calabi-Yau 3-folds 
and Aganagic-Vafa A-branes. 
In Section \ref{sec:openGWouter}, we introduce and compare
various GW invariants:
formal relative GW invariants of $(\hY,\hD)$, 
descendant GW invariants of $\hX$, 
open GW invariants of $(X,L)$,
and descendant GW invariants of $X$.
In Section \ref{sec:conj}, we review the geometry
of mirrors of toric Calabi-Yau 3-folds and Aganagic-Vafa
B-branes, and state the mirror conjecture on holomorphic
disks (Conjecture \ref{AVconj}). The proof of Conjecture
\ref{AVconj} is given in Section \ref{sec:proof}. We 
list explicit mirror formulae in Section \ref{sec:formula}.

\subsection*{Acknowledgments}
We thank C. Woodward, E. Zaslow, and J. Zhou for their comments.

\section{Geometry and Topology of Toric Calabi-Yau 3-Folds}
\label{sec:geometry}

\subsection{Toric varieties as geometric quotients}
We refer to \cite{Fu} for the theory of toric varieties. 
In this subsection, we consider smooth, possibly noncompact toric varieties of any dimension.

Let $N \cong \bZ^n$ be a free abelian group, and let $\tri$ be a \emph{fan} in $N$ (or in $N_\bR = N \otimes \bR\cong \bR^n$) 
of strongly convex rational polyhedral cones.  Let $X=X(\tri)$ be the toric variety associated to $\tri$. In this paper,
we assume that $X$ is smooth. We use the following notation:
\begin{itemize}
\item $\tri(d)$ is the set of $d$-dimensional cones in $\tri$.
\item Given $\si\in \tri(d)$, let $V(\si)$ denote the codimension $d$ orbit closure associated to $\si$.
\item Let $\tri(1)=\{\rho_1,\ldots, \rho_r\}$ be the set of 1-dimensional cones in $\tri$, and let $v_i\in N$ be 
the unique generator of the semigroup $\rho_i\cap N$, so that $\rho_i\cap N=\bZ_{\geq 0}v_i$.
\item Let $M=\Hom(M,\bZ)$ be the dual lattice of $N$.
\end{itemize}
There is a group homomorphism
$$
\phi: \tN:=\bigoplus_{i=1}^r\bZ \tv_i \cong \bZ^r \lra N\cong \bZ^n,\quad \tv_i\mapsto v_i.
$$
We assume that $\phi$ is surjective. 
Let 
$$
l^{(a)} = (l^{(a)}_1,\ldots, l^{(a)}_r), \quad a=1,\ldots,k, 
$$
be a $\bZ$-basis of $\Ker\phi\cong \bZ^k$, where $k=r-n$.
Let $\tM=\Hom(M;\bZ)$ be the dual lattice of $\tN$. Since
$\phi: \tN\to N$ is surjective,  the dual map $\phi^*:M\to \tM$ is injective. 
Applying $\Hom(-,\bC^*)$ to the following short exact sequence 
$$
0\to M\stackrel{\phi^*}{\lra} \tM \to A_{n-1}(X)\to 0,
$$
we obtain another short exact sequence
\begin{equation}\label{eqn:T}
1\to G\to \tbT \to \bT\to 1,
\end{equation}
where 
\begin{eqnarray*}
G &:=&\Hom(A_{n-1}(X),\bC^*)\cong (\bC^*)^k, \\
\tbT &=& \Hom(\tM,\bC^*)\cong (\bC^*)^r, \\
\bT &=& \Hom(M,\bC^*)\cong (\bC^*)^n.
\end{eqnarray*}
We also denote
$$\bL=\Ker \phi, \quad \bL^\vee=A_{n-1}(X).$$
Notice that $\bL^\vee \cong H^2(X;\bZ)$, and $\bL\cong H_2(X;\bZ)$.

The torus $\tbT\cong (\bC^*)^r$ acts on $\bC^r=\Spec\bC[X_1,\ldots, X_r]$.
Let $I(\tri) \subset \bC[X_1,\ldots, X_r]$ be
the ideal generated by $\{ \prod_{\rho_i\not \subset \si} X_i\mid \sigma\in \tri\}$,
and let  $Z(\tri)\subset \bC^r$ be the subvariety defined by $I(\tri)$.
Then $X$ can be identified with the geometric quotient:
\begin{equation}\label{eqn:geometric-quotient}
X = (\bC^r-Z(\tri))/G.
\end{equation}
where $G=(\bC^*)^k$ acts on $\bC^r$ by 
$$
(t_1,\ldots, t_k)\cdot (X_1,\ldots, X_r)
=(\prod_{a=1}^k t_a^{l_1^{(a)}}\cdot X_1,\ldots, 
\prod_{a=1}^k t_a^{l_r^{(a)}} \cdot X_r).
$$

The $\tbT$-divisor $\{X_i =0\}$ in $\bC^r-Z(\tri)$ descends to
a $\bT$-divisor $D_i$ in $X$. When $X$ is semi-projective, the quotient \eqref{eqn:geometric-quotient}
is also a quotient in the sense of GIT (geometric invariant theory) 
\cite[Section 2]{HS}.

\subsection{Toric varieties as symplectic quotients} 
When $X$ is a smooth semi-projective toric variety, we may describe K\"{a}hler structures on $X$ 
as a symplectic quotient. Let $G_\bR \cong U(1)^k$ be the maximal compact subgroup of
$G\cong (\bC^*)^k$. Then the dual $\fg_\bR^*$ of the Lie algebra $\fg_\bR$ of $G_\bR$
can be canonically identified with $A_{n-1}(X)\otimes \bR=H^{1,1}(X_\Si;\bR)$. 
Let $\tmu:\bC^r\to \fg_\bR^*\cong\bR^k$ be the moment map of the Hamiltonian $G_\bR$-action 
on $\bC^r$, equipped with the standard K\"{a}hler form
\begin{equation}\label{eqn:CrKahler}
\frac{\sqrt{-1}}{2} \sum_{i=1}^r dX_i\wedge d\overline{X}_i.
\end{equation}
Then 
$$
\tmu(X_1,\ldots, X_r) =
(\sum_{i=1}^r l^{(1)}_i|X_i|^2,\ldots, \sum_{i=1}^r l^{(k)}_i |X_i|^2).
$$
Let $(r_1,\ldots, r_k) \in H^{1,1}(X;\bR) \cong \bR^k$ be a K\"{a}hler class. Then
\begin{equation}\label{eqn:symplectic-quotient}
X = \tmu^{-1}(r_1,\ldots, r_k)/G_\bR.
\end{equation}
The standard K\"{a}hler form \eqref{eqn:CrKahler} on $\bC^r$ descends
to a K\"{a}hler form $\omega_{r_1,\ldots,r_k}$ on the quotient \eqref{eqn:symplectic-quotient}
in class $(r_1,\ldots,r_k)\in H^{1,1}(X;\bR)$.
The real numbers $r_1,\ldots,r_k$ are known as {\em K\"{a}hler parameters}
of $X$. Let $T_a = -r_a +\sqrt{-1}\theta_a$ be {\em complexified K\"{a}hler parameters} of $X$.

\subsection{The Calabi-Yau condition}
We say a toric variety $X$ is Calabi-Yau if the canonical divisor
$K_X = -D_1-\cdots -D_r$ is trivial. $X$ is Calabi-Yau  if and only if 
$$
\sum_{i=1}^r l_i^{(a)}=0,\quad a=1,\ldots,k 
$$

When $X$ is Calabi-Yau, we have a short exact sequence
\begin{equation}\label{eqn:Tzero}
1\to \bT'\to \bT\to \bC^*\to 1,
\end{equation}
where $\bT'\cong (\bC^*)^{n-1}$ is the subtorus of $\bT$
that acts trivially on $\Lambda^n T_X \cong \cO_X$.

\subsection{Aganagic-Vafa A-branes} \label{sec:A-brane}

Let $X =\tmu^{-1}(r_1,\ldots,r_k)/G_\bR$ be a smooth semi-projective toric 
Calabi-Yau 3-fold equipped with the K\"{a}hler form $\omega_{r_1,\ldots,r_k}$. Here
$\tmu^{-1}(r_1,\ldots, r_k) $ is defined by
$$
\sum_{i=1}^{k+3} l^{(a)}_i |X_i|^2 = r_a,\quad a=1,\ldots, k.
$$
Write $X_i=\rho_i e^{\sqrt{-1}\phi_i}$, where $\rho_i=|X_i|$. 
In \cite{AV}, Aganagic-Vafa introduced Lagrangian submanifolds $L$ of
$X_\Si$ defined by additional constraints:
$$
\sum_{i=1}^{k+3} \hat{l}_i^1|X_i|^2 =c_1,\quad
\sum_{i=1}^{k+3} \hat{l}_i^2 |X_i|^2 =c_2,\quad 
\sum_{i=1}^{k+3} \phi_i = \mathrm{const},
$$
where
$$
\hat{l}^\alpha_i\in \bZ,\quad \sum_{i=1}^{k+3} \hat{l}^\alpha_i =0, \quad \alpha=1,2.
$$
Such a Lagrangian submanifold is diffeomorphic to $S^1\times\bC$
and intersects a unique 1-dimensional  orbit closure $V(\tau_L)$, where $\tau_L\in \tri(2)$, 
along a circle. We say $L$ is an {\em outer} brane if 
$V(\tau_L) \cong \bC$, and we say $L$ is an {\em inner} brane if
$V(\tau_L) \cong \bP^1$.

\subsection{Framing} \label{sec:framing}
Let $\bT'_\bR\cong U(1)^2$ be the maximal compact subgroup of
$\bT'\cong (\bC^*)^2$. Then $\bT'_\bR$ preserves any Aganagic-Vafa A-brane
$L$. The $\bT'_\bR$-action on $L$ is given by 
$$
(t_1,t_2)\cdot (z,u) = (t_1z, t_2u)
$$
where $(t_1,t_2)\in \bT'_\bR= U(1)\times U(1)$, $(z,u)\in L\cong S^1\times \bC$.
Therefore $(X,L)$ is a $\bT'_\bR$-equivariant pair. The $\bT'_\bR$-action 
on $X$ is Hamiltonian with respect to the symplectic form $\omega_{r_1,\ldots,r_k}$. 
Let $\mu':X\to (\ft'_\bR)^*$ be the moment map of the $\bT'_\bR$-action, where
$(\ft'_\bR)^*\cong\bR^2$ is the dual of the Lie algebra $\ft'_\bR$ of $\bT'_\bR$.

We define the 1-skeleton of $X$ to be
$$
X^1=\bigcup_{\tau\in \tri(2)}V(\tau).
$$
Then $\mu'(X^1)$ is an immersed trivalent graph. 
$\mu'(V(\tau))$ is a line segment if $V(\tau)\cong\bP^1$, and is a ray
if $V(\tau)\cong \bC$.  
Suppose that $\si$ is a 3-dimensional cone in $\tri$ and $\tau_L \subset \si$. 
Then $V(\si) \subset V(\tau_L)$, where $V(\si)$ is the $\bT$-fixed point associated to $\si$. Let 
$$
\wu_L = (c_1)_{\bT'}(T_{V(\si)} V(\tau_L)) \in H^2_{\bT'}(V(\si);\bZ) \cong H^2(B\bT';\bZ)\subset H^2(B\bT';\bR)
\cong (\ft'_\bR)^*. 
$$
Then $\wu_L$ is tangent to $\mu(V(\tau_L))$.  Let $\tau'_L \in \tri(2)$ be the 2-cone
such that $\tau'_L\subset\si$ and  $\mu'(V(\tau'_L))$ is the first edge encountered in the counterclockwise direction from $\mu'(V(\tau_L))$. 
Define
$$
\wv_L= (c_1)_{\bT'}(T_{V(\si)}V(\tau'_L)) \in H^2_{\bT'}(V(\si);\bZ) \subset (\ft'_\bR)^*.
$$

Let $\rho_{i_1}\subset \sigma$ be the $1$-cone that is not in $\tau_L$, and let $\rho_{i_2}\subset \sigma$ be the $1$-cone that is not in $\tau_L'$. We use $\rho_{i_3}$ to denote the $1$-cone in $\sigma$ other than $\rho_{i_1},\rho_{i_2}$. If $V(\tau_L)\cong \bP^1$, then there is another $3$-cone $\sigma'\in \tri(3)$ other 
than $\sigma$ of which $\tau_L$ is a face. Let $\rho_{i_4}$ be the $1$-cone in $\sigma'$ but not in $\tau_L$.

\begin{figure}[h]
\psfrag{1}{\small $D_{i_1}$} 
\psfrag{2}{\small $D_{i_3}$}
\psfrag{3}{\small $D_{i_2}$} 
\psfrag{4}{\small $D_{i_4}$}
\psfrag{p}{\footnotesize $V(\sigma)$}
\psfrag{q}{\footnotesize $V(\sigma')$}
\includegraphics[scale=0.6]{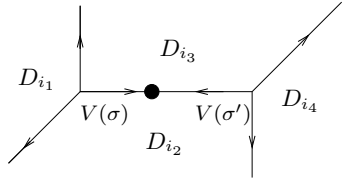}
\caption{images of $D_{i_1},D_{i_2},D_{i_3},D_{i_4}$, $V(\si)$, $V(\si')$ under $\mu'$.}
\label{fig:inner}
\end{figure} 

\begin{figure}[h] 
\psfrag{u}{\small $\wu_L$}
\psfrag{v}{\small $\wv_L$}
\psfrag{-u-v}{\small $-\wu_L -\wv_L$}
\psfrag{v-fu}{\small $\wf=\wv_L- f\wu_L$}
\psfrag{fu-v}{\small $-\wf=f\wu_L - \wv_L$}
\psfrag{mu(L)}{\small $\mu'(L)$}
\includegraphics[scale=0.6]{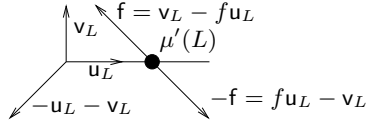}
\caption{framing of an A-brane}
\label{fig:L}
\end{figure}

Let $L$ be an Aganagic-Vafa A-brane. Then $\mu'(L)$ is a point. 
A framing of $L$ is a choice of a vector $\wf \in H^2(B\bT';\bZ)$ such that
$\wu_L \wedge \wf = \wu_L\wedge \wv_L$, i.e., $\wf= \wv_L - f \wu_L$ for some $f\in \bZ$.
We sometimes call the integer $f$ the a framing of $L$. See Figure \ref{fig:L}.

We may view $\wf$ as an element in $\Hom(\bT',\bC^*)$. Let
$\bT_{L,f}\cong \bC^*$ be the kernel of $\wf: \bT'\to \bC^*$.

\subsection{Topology} \label{sec:topology}

Suppose that $p_1,\ldots, p_k$ is a $\bZ$-basis of $H^2(X;\bZ)\cong \bL^\vee$ dual to $l^{(1)},\dots, l^{(k)}\in \bL$. We choose appropriate $\{l^{(a)}\}_{a=1}^k$ and $\{p_a\}_{a=1}^k$ such that the 
class of $\omega_{r_1,\ldots,r_k}$ is  $\sum_{a=1}^k r_a p_a\in H^2(X_\Si;\bR)=H^{1,1}(X_\Si;\bR)$. 
Let $Q_a= e^{T_a}$, which are formal K\"{a}hler parameters. Given $\beta\in H_2(X;\bZ)$, define
$$
Q^\beta:= e^{\langle \beta, \sum_{a=1}^k T_a p_a \rangle} =\prod_{a=1}^k Q_a^{d_a},
$$
where $d_a=\langle \beta, p_a\rangle$, $a=1,\ldots,k$.
For $i=1,\ldots,r$, we have
$$
D_i^*:= c_1(\cO_X(D_i))=\sum_{a=1}^k l^{(a)}_i p_a, \quad
\langle D_i^*,  \beta\rangle =\sum_{a=1}^k  d_a l_i^{(a)}.
$$

Let $\tri(2)_c=\{ \tau\in \tri(2)\mid V(\tau)\cong \bP^1\}$.
The inclusion $X^1\hookrightarrow X$ induces a surjective
group homomorphism 
$$
H_2(X^1;\bZ)=\bigoplus_{\tau\in \tri(2)_c} V(\tau) \to H_2(X;\bZ).
$$

Let 
$$
L=L_1\cup \cdots \cup L_s
$$
be a disjoint union of $s$ Aganagic-Vafa A-branes
$L_1,\ldots, L_s$ in $X$.  
We have the following long exact sequence of relative homology groups 
of the pair $(X, L)$:
\begin{equation} \label{eqn:XL}
\cdots \to H_2(L;\bZ) \to H_2(X;\bZ)\to H_2(X,L;\bZ)\stackrel{\partial}{\to} H_1(L;\bZ)
\to H_1(X;\bZ)\to \cdots
\end{equation}
where $H_2(L;\bZ)=0$ and $H_1(X;\bZ)=0$.

For $i=1,\ldots, s$, let $V(\tau_i)$ be the unique 1-dimensional orbit
closure that intersects $L_i$.  We assume that $\tau_1,\ldots, \tau_s$
are distinct. We choose a $\bT$-fixed point $x_i \in V(\tau_i)$.
When $L_i$ is an outer brane, $V(\tau_i) \cong \bC$, and there is
a unique such point; when $L_i$ is an inner brane, there
are exactly two such points, $x_i$ and $x_i^-$. Then $V(\tau_i)-L_i\cap V(\tau_i)$ has
two connected components, one of which is a holomorphic
disk $B_i$ containing $x_i$. We orient
$B_i$ by the holomorphic structure. Then $B_i$ represents a
relative homology class $b_i \in H_2(X,L;\bZ)$. 
Let $\gamma_i\in H_1(L_i;\bZ)$ be the image of
$b_i$ under the map $\partial$ in \eqref{eqn:XL}. 
Then 
$$
H_1(L;\bZ)=\bigoplus_{i=1}^s H_1(L_i;\bZ) =\bigoplus_{i=1}^s \bZ \gamma_i.
$$
Any element $\beta' \in H_2(X,L;\bZ)$  is of the form 
$$
\beta' =\beta+ \sum_{i=1}^s d_i b_i,
$$
where $d_i \in \bZ$ and $\beta\in H_2(X;\bZ)$.
	
Let $\cA$ be the collection of \emph{anti-cones}
$$
\cA=\{J\subset \{1,\dots, r\}|\text{$\sum_{j\notin J} \bR_{\ge 0} v_j$ is a cone in $\Delta$}\}.
$$
The effective cone $\bL_\eff$ is
$$
\bL_\eff=\{\beta\in \bL| \{i\in\{1,\dots, r\} |\langle D_i^*, \beta \rangle\in \bZ_{\ge 0}\}\in \cA\}.
$$
Usually the cone $\bL_\eff \otimes_\bZ \bR$ in the real vector space $\bL\otimes_\bZ \bR$ is called the Mori cone. One says $\beta\ge 0$ if $\beta \in \bL_\eff$ and $\beta>0$ if $\beta$ is also non-zero. The closure of the
K\"ahler cone in $\bL^\vee\otimes_\bZ \bR$ is 
$$
C_X=\bigcap_{I\in \cA} (\sum_{i\in I} \bR_{\ge 0} v_i).
$$
It is dual to the Mori cone.

\section{Gromov-Witten Invariants }
\label{sec:openGWouter}

Let $X$ be a toric Calabi-Yau 3-fold defined by a fan $\tri$. Let 
$$
(L_1,f_1),\ldots, (L_s,f_s)
$$ 
be framed Aganagic-Vafa A-branes (see Section \ref{sec:A-brane}), and
let $\tau_1,\ldots,\tau_s\in \tri(2)$ be defined as in Section \ref{sec:topology}.
Let $\si_1,\ldots,\si_s \in \tri(3)$ be the top dimensional cones
corresponding to $\bT$ fixed points $x_1,\ldots, x_s$, respectively, where
$x_i$ are defined as in Section \ref{sec:topology}. Then
$$
\tau_i\subset \si_i,\quad
V(\si_i)=\{x_i\},\quad i=1,\ldots,s.
$$ 
Define 
$$
w^i_1 =\wu_{L_i},\quad w^i_2 =\wv_{L_i},\quad w^i_3 = -w^i_1-w^i_2,
\quad \wf_i = w^i_2 -f_i w^i_1.
$$

In Section \ref{sec:YD}--\ref{sec:multiple}, we assume $L_1,\ldots, L_s$ are outer branes.
In Section \ref{sec:single} --\ref{sec:generating}, $L=L_1$ is a single outer or inner brane.

\subsection{The relative FTCY 3-fold $(\hY,\hD)$} \label{sec:YD}
\begin{figure}[h]
\psfrag{w1}{\small $\ww^i_1$}
\psfrag{w2}{\small $\ww^i_2$}
\psfrag{w3=-w1-w2}{\small $\ww^i_3=-\ww^i_1-\ww^i_2$}
\psfrag{f=w2-fw1}{\small $\wf_i =\ww^i_2 - f_i \ww^i_1$}
\psfrag{-f=fw1-w2}{\small $-\wf_i= f_i \ww^i_1-\ww^i_2$}
\psfrag{e}{\small $\bar{e}_i$}
\psfrag{ep}{\small $\bar{e}_i'$}
\psfrag{v}{\small $v_i$}
\psfrag{vp}{\small $v'_i$}
\psfrag{GaX}{\small $\Gamma_X$}
\psfrag{GaXL}{\small $\Gamma_{X,(L_1,f_1),\ldots,(L_k,f_k)}$}
\includegraphics[scale=0.6]{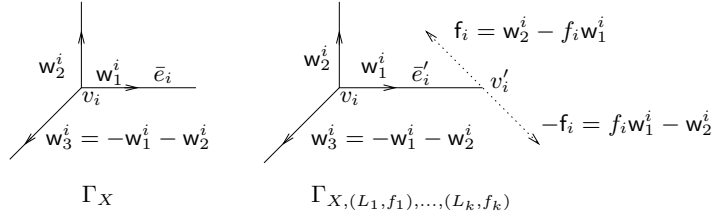}
\label{fig:FTCY}
\caption{The FTCY graphs $\Gamma_X$ and $\Gamma_{X,(L_1,f_1),\ldots, (L_s,f_s)}$. The trivalent vertex
$v_i$ corresponds to the cone $\si_i\in \tri(3)$.}
\end{figure}

We refer to \cite[Section 3]{LLLZ} for the definition
of formal toric Calabi-Yau (FTCY) graphs and
construction of relative FTCY 3-folds.
Let $\Gamma_X$ be the formal toric Calabi-Yau (FTCY) graph associated to 
the smooth toric 3-fold $X$ (see \cite[Section 3.1]{LLLZ}). 
We define a FTCY graph 
$$
\Gamma_{X,(L_1,f_1),\ldots, (L_s,f_s)}
$$ 
by replacing the noncompact edge $\bae_i$ in $\Gamma_X$ associated to
$\tau_i$ by a compact edge $\bae_i'$  with framing $\wf_i$.\footnote{In \cite{LLLZ}
all the noncompact edges are replaced by a compact edge ending at an univalent
vertex. Here we only compactify $\bae_1,\ldots, \bae_s$.}
Figure \ref{fig:FTCY} shows the construction near the edge $\bae_i$. 

Let $(\hY,\hD)$ be the relative FTCY 3-fold associated to
the FTCY graph $\Gamma= \Gamma_{X,(L_1,f_1),\ldots, (L_s,f_s)}$. Then
$\hD$ is the disjoint union of its connected components
$\hD^1,\ldots, \hD^s$ which correspond to the
$s$ univalent vertices $v'_1,\ldots,v'_s$ in $\Gamma$.
$\bT=(\bC^*)^3$ and $\bT'=(\bC^*)^2$ acts on $(\hY,\hD)$.
Let $E_c(\Gamma_X)$ and $E_c(\Gamma)$ denote the
set of compact edges in $\Gamma_X$ and $\Gamma$, respectively. 
Then  $E_c(\Gamma)=E_c(\Gamma_X)\cup \{ \bae'_1,\ldots, \bae'_s\}$.
Each compact edge $\bae$ corresponds to a torus invariant
projective line $C^{\bae}$ in $X$ or in $\hY$. We have
\begin{eqnarray*}
H_2(X^1;\bZ) &=& \bigoplus_{\bae\in E_c (\Gamma_X)} \bZ [C^{\bae}],\\
H_2(\hY;\bZ) &=& \bigoplus_{\bae\in E_c(\Gamma)} \bZ[C^{\bae}] = H_2(X^1;\bZ)\oplus \bigoplus_{i=1}^s \bZ[C^{\bae_i}],
\end{eqnarray*}
where $C^{\bae_i}=V(\tau_i)$. There is a surjective group homomorphism
$$
\pi: H_2(\hY;\bZ) = H_2(X^1;\bZ)\oplus \bigoplus_{i=1}^s \bZ[C^{\bae_i}] \longrightarrow
H_2(X,L;\bZ)=H_2(X;\bZ)\oplus \bigoplus_{i=1}^s \bZ b_i
$$
given by 
$$
\beta+\sum_{i=1}^s d_i[C^{\bae_i}] \mapsto
i_*\beta + \sum_{i=1}^s d_i b_i,
$$
where $\beta\in H_2(X^1;\bZ)$, and $i_*:H_2(X^1;\bZ)\to H_2(X;\bZ)$ is induced by the inclusion $i:X^1\hookrightarrow X$.

\subsection{Formal relative GW invariants of $(\hY,\hD)$}  \label{sec:GWYD}

We refer to \cite[Section 4]{LLLZ} for the precise definitions
of relative stable morphisms to a relative FTCY 3-fold, and 
formal relative GW invariants. Here we
briefly recall the definitions for our particular setting.

An {\em effective class} of $\Gamma$ is a pair $(\vd',\vmu)$, where
\begin{itemize}
\item $\vd': E_c(\Gamma)\to \bZ_{\geq 0}$ is a map from the set of compact edges in $\Gamma$ to the
set of nonnegative integers,
\item $\vmu: \{ v'_1,\ldots, v'_s\} \to \cP$ is a map from the
set of univalent vertices in $\Gamma$ to set of partitions,
\item $\vmu(v'_i)$ is a partition of $\vd'(\bae'_i)$.
\end{itemize}
$\vmu$ can be viewed as an $s$-uple of partitions: $\vmu=(\mu^1,\ldots,\mu^s)$, where $\mu^i :=\vmu(v'_i)$.

Define
$$
\ell(\vmu)=\sum_{i=1}^s \ell(\mu^i),\quad \Aut(\vmu)=\prod_{i=1}^s\Aut(\mu^i).
$$
We also use $\vd'$ to denote the following effective curve class:
$$
\sum_{\bae\in \Gamma} \vd'(\bae)[C^{\bae}] \in H_2(\hY;\bZ).
$$
A genus $g$, class $(\vd',\vmu)$ relative stable morphism to $(\hY,\hD)$ is
an $\ell(\vmu)$-pointed, genus $g$, degree $\vd'$ stable morphism $u:(C,\bq)\to \hY$, where
$$
\bq=\{ q^i_j\mid 1\leq i\leq s,1\leq j\leq \ell(\mu^i)\}
$$
are distinct smooth points on $C$, and
$$
u^{-1}(\hD^i)=\sum_{i=1} \mu^i_j q^i_j
$$
as Cartier divisors. Let $\cM_{g,\vd',\vmu}(\hY,\hD)$ be the moduli space of genus $g$, class $(\vd',\vmu)$
relative stable morphisms to $(\hY,\hD)$. 
The stable compactification $\Mbar_{g,\vd',\vmu}(\hY,\hD)$
of $\cM_{g,\vd',\vmu}(\hY,\hD)$ contains relative morphisms
to expanded targets $(\hY_\bm,\hD_\bm)$. 
Then $\Mbar_{g,\vd',\vmu}(\hY,\hD)$ is a formal Deligne-Mumford stacks equipped with 
a perfect obstruction theory of virtual dimension $\ell(\vmu)$.  Moreover,
$\bT'=(\bC^*)^2$ acts $\Mbar_{g,\vd',\vmu}(\hY,\hD)$, and the perfect obstruction theory are $\bT'$-equivariant.
The genus $g$, class $(\vd',\vmu)$ formal relative GW invariants of $(\hY,\hD)$ is defined by 
\begin{equation} \label{eqn:Ngdmu}
F^{\hY,\hD}_{g,\vd',\vmu} = \frac{1}{ |\Aut(\vmu)| }\int_{[\Mbar_{g,\vd',\vmu}(\hY,\hD)^{\bT'}]^\vir }
\frac{\prod_{i=1}^s\prod_{j=1}^{\ell(\mu^i)} (\ev^i_j)^*(c_1)_{\bT'}(\cO_{\hD_i}(\hL_i))}{e_{\bT'}(N^\vir)}, 
\end{equation}
where $\Mbar_{g,\vd',\vmu}(\hY,\hD)^{\bT'}$ is the $\bT'$ fixed points substack in 
$\Mbar_{g,\vd',\vmu}(\hY,\hD)$,  $\ev^i_j$ is the evaluation at
the marked point $q^i_j$, $\hL^i$ is a $\bT'$-equivariant divisor in $\hD^i$, 
and $N^\vir$ is the virtual normal bundle
of $\Mbar_{g,\vd',\vmu}(\hY,\hD)^{\bT'}$ in $\Mbar_{g,\vd',\vmu}(\hY,\hD)$. 
Let $\{ \wu, \wv\}$ be a $\bZ$-basis of $H^2(B\bT';\bZ)$, so that $H^*(B\bT';\bQ)=\bQ[\wu,\wv]$.
A priori the right hand side of \eqref{eqn:Ngdmu} is a rational function in $\wu,\wv$, homogeneous
of degree $0$, but it is indeed a rational number independent of $\wu,\wv$. 

Introduce variables $\bp=\{ p^i_j\mid 1\leq i\leq s, j\in \bZ_{>0}\}$ and $Q=\{Q_{\bae}\mid \bae \in E_c(\Ga)\}$. 
We define generating functions of formal relative GW invariants of $(\hY,\hD)$:
\begin{equation}
F^{\hY,\hD}_{\vd',\vmu}(\lambda) = \sum_g N^{\hY,\hD}_{\vd',\vmu} \lambda^{2g-2+\ell(\vmu)},
\end{equation}
\begin{equation}
F^{\hY,\hD}(\lambda,Q,\bp) = \sum_{\tiny \begin{array}{c} (\vd',\vmu)\\ \vmu\neq (\emptyset,\ldots,\emptyset)\end{array}} 
F^{\hY, \hD}_{\vd',\vmu}(\lambda)Q^{\vd'} \bp_{\vmu}
\end{equation}
where
$\displaystyle{Q^{\vd'}= \prod_{\bae\in E_c(\Gamma)} Q_{\bae}^{\vd'(\bae)} }$, 
$\displaystyle{\bp_{\vmu}= \prod_{i=1}^s \prod_{j=1}^{\ell(\mu^i)} p^i_{\mu^i_j} }$.
Define generating functions $F^{\bu\hY\hD}_{\vd',\vmu}$ and $Z^{\hY,\hD}$ by
\begin{equation}
Z^{\hY,\hD}(\lambda,Q,\bp) = \exp\left(F^{\hY,\hD}(\lambda,q,\bp)\right)= 
\sum_{(\vd',\vmu)} F^{\bullet \hY, \hD}_{\vd',\vmu}(\lambda)Q^{\vd'} \bp_{\vmu}.
\end{equation}

\subsection{Descendant GW invariants of $\hX$} \label{sec:GWhX} 
Let $\hX$ be the formal completion of $X$ along the 1-skeleton $X^1$. Then $\hX=\hY-\hD$
is a FTCY 3-fold, and 
$$
H_2(\hX;\bZ)=H_2(X^1;\bZ)=\bigoplus_{\bae\in E_c(\Gamma_X)} \bZ[C^{\bae}].
$$
The surjective group homomorphism $\pi:H_2(\hY;\bZ)\to H_2(X,L;\bZ)$ restricts
to a surjective group homomorphism $H_2(\hX;\bZ)\to H_2(X;\bZ)$.

An effective class of $\Gamma_X$ is  a function $\vd:E_c(\Gamma_X)\to \bZ_{\geq 0}$; we also
use $\vd$ to denote the effective curve class
$$
\sum_{\bae\in E_c(\Gamma_X)} \vd(\bae) [C^{\bae}] \in H_2(\hX;\bZ).
$$
Let $\Mbar_{g,n}(\hX,\vd)$ be the moduli space of
$n$-pointed, genus $g$, degree $\vd$ stable
morphisms to $\hX$.  Define
\begin{equation}\label{eqn:Ggd}
\begin{aligned}
G^{\hX}_{g,\vd, \vmu} (\wu,\wv) =& \frac{(-\sqrt{-1})^{\ell(\vmu)}}{ |\Aut(\vmu)| }
\prod_{i=1}^s \prod_{j=1}^{\ell(\mu^i)}
\frac{\prod_{m=1}^{\mu^i_j-1}(\mu^i_j \ww^i_2 +  m \ww^i_1) }{\mu^i_j\cdot \mu^i_j! (\ww^i_1)^{\mu^i_j-1}}. \\
& \int_{[\Mbar_{g,\ell(\vmu)}(\hX,\vd)^{\bT'}]^{\vir}}\frac{1}{e_{\bT'}(N^\vir)}
\prod_{i=1}^s\prod_{j=1}^{\ell(\mu^i)}\frac{ (\ev^i_j)^*\phi_i}{\frac{\ww_1^i}{\mu^i_j}(\frac{\ww_1^i}{\mu^i_j}-\psi^i_j)}
\end{aligned}
\end{equation}
where $\phi_i\in H_{\bT'}(\hX;\bZ)$ is the $\bT'$-equivariant
Poincar\'{e} dual of the torus fixed point $x_i\in V(\tau_i)$, and
$$
\ev^i_j: \Mbar_{g,\ell(\vmu)}(\hX,\vd)\to \hX, \quad i=1,\ldots, s,\quad j=1,\ldots, \ell(\mu^i),
$$
are evaluations at the marked points.

Define generating functions
\begin{eqnarray*}
G^{\hX}_{\vd,\vmu}(\lambda;\wu,\wv)&=& \sum_g \lambda^{2g-2+\ell(\vmu)}G^{\hX}_{\vd,\vmu}(\wu,\wv)\\
G^{\hX}(\lambda,Q,\bp;\wu,\wv)&=& 
\sum_{\tiny \begin{array}{c} (\vd,\vmu)\\ \vmu\neq(\emptyset,\ldots,\emptyset) \end{array}} 
G^{\hX}_{\vd,\vmu}(\lambda;\wu,\wv)Q^{\vd} \bp_{\vmu}\\
G^{\bullet\hX}(\lambda,Q,\bp;\wu,\wv)&=& \exp\left(G^{\hX}(\lambda,Q,\bp;\wu,\wv)\right)
= \sum_{(\vd,\vmu)} G^{\bullet\hX}_{\vd,\vmu}(\lambda;\wu,\wv)Q^{\vd} \bp_{\vmu} 
\end{eqnarray*}

Given two partitions $\mu^+, \mu^-$ of $d$, define
$$
\Phi_{\mu^+,\mu^-}^\bu(\lambda)=\sum_{\nu\vdash d} e^{\kappa_\nu \lambda/2}
\frac{\chi_\nu(\mu^+)}{z_{\mu^+}} \frac{\chi_\nu(\mu^-)}{z_{\mu^-}}
$$
where 
$\displaystyle{\kappa_\mu =\sum_{j=1}^{\ell(\mu)} \mu_j(\mu_j-2j+1) }$,
$\displaystyle{ z_\mu= \Aut(\mu)\prod_{j=1}^{\ell(\mu)}{\mu_j} }$,
and $\chi_\mu(\nu)$ denotes the value of the irreducible character of $S_d$ associated to $\nu\vdash d$
at the conjugacy class of $S_d$ associated to $\mu\vdash d$.

\begin{proposition} \label{pro:YDX}
\begin{equation}
\begin{aligned}
 & \sqrt{-1}^{\ell(\vmu)} (-1)^{\sum_{i=1}^s f_i(|\mu^i|-1)} F^{\bu\hY,\hD}_{\vd', \vmu} (\lambda)\\ 
=& \sum_{|\nu^i|=|\mu^i|}  
G^{\bu\hX}_{\vd, \vnu} (\lambda;\wu,\wv)
\prod_{i=1}^s z_{\nu^i}\Phi_{\nu^i, \mu^i}^\bu (\sqrt{-1}\frac{\wf_i}{\ww^i_1})\\
=& \sum_{|\nu^i|=|\mu^i|}  G^{\bu\hX}_{\vd, \vnu} (\lambda;\wu,\wv)
\prod_{i=1}^s z_{\nu^i}\Phi_{\nu^i, \mu^i}^\bu (\sqrt{-1}(f_i-\frac{\ww^i_2}{\ww^i_1})),
\end{aligned}
\end{equation}
where $\vd:E_c(\Gamma_X)\to \bZ_{\geq 0}$ is the restriction
of $\vd':E_c(\Gamma)\to\bZ_{\geq 0}$.
\end{proposition}
\begin{proof} We use the notation in \cite[Section 7]{LLLZ}, and set $\nu^i=\nu^{\bae_i}$.
By \cite[Proposition 7.10]{LLLZ},
\begin{eqnarray*}
F^{\bu\hY,\hD}_{\vd', \vmu} (\lambda)&=&\sum_{|\nu^{\bae}|=\vd'(\bae)}
\prod_{\bae\in E(\Ga)} (-1)^{n^e \vd'(\bae)} z_{\nu^{\bae}}
\prod_{v\in V_3(\Ga)}\sqrt{-1}^{\ell(\vnu^v)} G^\bu_{\vnu^v}(\lambda;\mathbf{w}_v) \\
&& \quad \cdot \prod_{i=1}^s (-1)^{|\mu^i|}(-\sqrt{-1})^{\ell(\nu^i) +\ell(\mu^i)} \Phi^\bu_{\nu^i,\mu^i}(\sqrt{-1}\frac{\wf_i}{\ww^i_1} \lambda)\\
& =&\sum_{|\nu^i|=|\mu^i|} \left(\sum_{|\nu^{\bae}|=\vd(\bae)} \prod_{\bae\in E(\Ga_X)} (-1)^{n^e\vd(\bae)} z_{\nu^{\bae}}
\prod_{v\in V_3(\Ga)}\sqrt{-1}^{\ell(\vnu^v)} G^\bu_{\vnu^v}(\lambda;\mathbf{w}_v)\right) \\
&& \prod_{i=1}^s (-1)^{(f_i-1)|\mu^i|}(-\sqrt{-1})^{\ell(\nu^i) +\ell(\mu^i)}z_{\nu^i} \Phi^\bu_{\nu^i,\mu^i}(\sqrt{-1}\frac{\wf_i}{\ww^i_1} \lambda)
\end{eqnarray*}
where 
$$
\sum_{|\nu^{\bae}|=\vd(\bae)} \prod_{\bae\in E(\Ga_X)} (-1)^{n^e}d^{\bae} z_{\nu^{\bae}}
\prod_{v\in V_3(\Ga)}\sqrt{-1}^{\ell(\vnu^v)} G^\bu_{\vnu^v}(\lambda;\mathbf{w}_v)
=\sqrt{-1}^{\ell(\vnu)}G^{\bu\hX}_{\vd,\vnu}(\lambda;\wu,\wv)
$$
by localization computations similar to those in \cite{LLLZ}. Therefore,
$$
F^{\bu\hY,\hD}_{\vd', \vmu} (\lambda)= (-\sqrt{-1})^{\ell(\vmu)} (-1)^{\sum_{i=1}^s(f_i-1)|\mu^i|}
\sum_{|\nu^i|=|\mu^i|}G^{\bu\hX}_{\vd,\vnu}(\lambda;\wu,\wv) z_{\nu^i}\Phi^\bu_{\nu^i,\mu^i}(\sqrt{-1}\frac{\wf_i}{\ww^i_1} \lambda)
$$

\end{proof}

\subsection{Open GW invariants of $(X,L)$: multiple outer branes}\label{sec:multiple}
The open GW invariants of $X$ relative to the framed outer branes
$(L_1,f_1),\ldots, (L_k,f_k)$ are defined to be
\begin{equation}\label{eqn:XLYD}
N^{X,L}_{g,\beta',\mu^1,\ldots,\mu^s}(f_1,\ldots,f_s) = (-1)^{\sum_{i=1}^s(|\mu^i|-\ell(\mu^i))}
\sum_{\pi(\vd')=\beta'} F^{\hY,\hD}_{g,\vd',\vmu} \in \bQ,
\end{equation}
where $\vmu(v_i')=\mu^i$ for $i=1,\ldots,s$, and the sign $(-1)^{\sum_{i=1}^s(|\mu^i|-\ell(\mu^i) ) }$
is included for convenience of comparison with the predictions from mirror symmetry.

Given $\beta\in H_2(X;\bZ)$ and $\vmu=(\mu^1,\ldots, \mu^s)$, define
\begin{equation}\label{eqn:Ggbeta}
\begin{aligned}
G^X_{g,\beta, \vmu} (\wu,\wv) =& \frac{(-\sqrt{-1})^{\ell(\vmu)}}{ |\Aut(\vmu)| }
\prod_{i=1}^s \prod_{j=1}^{\ell(\mu^i)}
\frac{ \prod_{m=1}^{\mu^i_j-1}(\mu^i_j \ww_2^i +  m \ww_1^i) }{\mu^i_j \cdot \mu^i_j! (\ww_1^i)^{\mu^i_j-1} } \\
&\cdot  \int_{[\Mbar_{g,\ell(\vmu)}(X,\beta)^{\bT'}]^{\vir}}\frac{1}{e_{\bT'}(N^\vir)}
\prod_{i=1}^s\prod_{j=1}^{\ell(\mu^i)}\frac{(\ev^i_j)^*\phi_i}{\frac{\ww_1^i}{\mu^i_j}(\frac{\ww_1^i}{\mu^i_j}-\psi^i_j)}
\end{aligned}
\end{equation}
Then
$$
G^X_{g,\beta,\vmu}=\sum_{\pi(\vd)= \beta} G^{\hX}_{g,\vd,\vmu}
$$
In particular, when $X=\bC^3$, $s=3$, we have
$$
\ww^1_1=\ww^3_2 =\wu,\quad
\ww^1_2 = \ww^2_1 =\wv,\quad
\ww^2_2 = \ww^3_1 = -\wu-\wv,
$$
$$
\Mbar_{g,\ell(\vmu)}(\bC^3,0)^{\bT'}=\Mbar_{g,\ell(\vmu)},
$$
$$
\frac{1}{e_{\bT'}(N^\vir)} =  \prod_{i=1}^3 \frac{\Lambda_g^\vee(\ww_i)}{\ww_i} \in H^*(\Mbar_{g,\ell(\vmu)};\bQ)\otimes \bQ(\wu,\wv), 
$$
where $\Lambda_g^\vee(u)=u^g -\lambda_1 u^{g-1} +\cdots +(-1)^g\lambda_g$. Define
$\ww_1=\ww_4=\wu$, $\ww_2 = \wv$, and $\ww_3 =-\wu-\wv$. Then
\begin{eqnarray*}
G^{\bC^3}_{g,0,\vmu} (\wu,\wv) &=& \frac{(-\sqrt{-1})^{\ell(\vmu)}}{ |\Aut(\vmu)| }
\prod_{i=1}^3 \prod_{j=1}^{\ell(\mu^i)}
\frac{ \prod_{m=1}^{\mu^i_j-1}(\mu^i_j \ww_{i+1} +  m \ww_i) }{\mu^i_j \cdot \mu^i_j! \ww_i^{\mu^i_j-1} } \\
&&\cdot\int_{\Mbar_{g,\ell(\vmu)}}\prod_{i=1}^3\frac{\Lambda_g^\vee(\ww_i)\ww_i^{\ell(\vmu)-1} }
{\prod_{m=1}^{\mu^i_j  }\frac{\ww_i}{\mu^i_j}(\frac{\ww_i}{\mu^i_j}-\psi^i_j) } 
\end{eqnarray*}
which is exactly the three-partition Hodge integral defined in \cite[Section 2.3]{LLLZ}.

Given $\beta'\in H_2(X,L;\bZ)$ and an $s$-uple of partitions
$\vmu= (\mu^1,\ldots,\mu^s)$, define 
$$
F^{\bu X,L}_{\beta',\vmu}(\lambda;f_1,\ldots,f_s) 
=(-1)^{\sum_{i=1}^s(|\mu^i|-\ell(\mu^i)) } \sum_{\pi(\vd')=\beta'}F^{\bu \hY,\hD}_{\vd',\vmu}(\lambda).
$$

The following proposition follows from Proposition \ref{pro:YDX} and definitions.
\begin{proposition}[multiple framed outer branes] \label{pro:XLX}
Let $\beta\in H_2(X;\bZ)$ be an effective curve class, and let
$\beta'=\beta+\sum_{i=1}^s|\mu^i|b_i \in H_2(X,L;\bZ)$. Then
\begin{equation} \label{eqn:tF-GX}
\begin{aligned}
& (-\sqrt{-1})^{\ell(\vmu)} (-1)^{\sum_{i=1}^s f_i |\mu^i|}F^{\bu X,L}_{\beta', \vmu} (\lambda;f_1,\ldots,f_s) \\
=& \sum_{|\nu^i|=|\mu^i|}  
G^{\bu X}_{\beta, \vnu} (\lambda;\wu,\wv)
\prod_{i=1}^s z_{\nu^i}\Phi_{\nu^i, \mu^i}^\bu (\sqrt{-1}\frac{\wf_i}{\ww^i_1})\\
=& \sum_{|\nu^i|=|\nu^i|}  G^{\bu X}_{\beta, \vnu} (\lambda;\wu,\wv)
\prod_{i=1}^s z_{\nu^i}\Phi_{\nu^i, \mu^i}^\bu (\sqrt{-1}(f_i-\frac{\ww^i_2}{\ww^i_1})).
\end{aligned}
\end{equation}
\end{proposition}

In particular, when $s=1$, $\mu^1=\mu$, $f_1=f$, $b_1=b$, we have
$$
F^{\bu X,L}_{\beta+|\mu|b, \mu} (\lambda;f) =\sqrt{-1}^{\ell(\mu)}(-1)^{f|\mu|}
\sum_{|\nu|=|\mu|}  G^{\bu X}_{\beta, \nu} (\lambda;\wu,\wv)
z_{\nu}\Phi_{\nu, \mu}^\bu (\sqrt{-1}(f-\frac{\wv_L}{\wu_L})).
$$

Let $\bT_{L,f}\cong \bC^*$ be the subtorus of $\bT'$ defined
in Section \ref{sec:framing}. It is the kernel of
the character $\wv_L-f\wu_L\in \Hom(\bT',\bC^*)$. The inclusion
$\bT_{L,f}\hookrightarrow \bT'$ induces a surjective ring homomorphism 
$$
H^*(B\bT';\bQ)=\bQ[\wu,\wv]= \bQ[\wu_L,\wv_L] \to  H^*(\bT_{L,f};\bZ)=\bQ[\wu_L]
$$
given by $\wu_L\mapsto \wu_L$ and $\wv_L\mapsto f\wu_L$.

Recall that
$$
\Phi_{\nu,\mu}^{\bu}(0)=\frac{\delta_{\nu,\mu}}{z_\nu},
$$
so 
$$
F^{\bu X,L}_{\beta +|\mu|b , \mu} (\lambda;f) =\sqrt{-1}^{\ell(\mu)}(-1)^{f|\mu|}
\sum_{|\nu|=|\mu|}  G^{\bu X}_{\beta, \nu} (\lambda;\wu,\wv)\Bigr|_{\wv_L=f\wu_L}.
$$
Therefore,
$$
N_{g,\beta +|\mu|b ,\mu}^{X,L}(f)= \sqrt{-1}^{\ell(\mu)}(-1)^{f|\mu|}G^X_{\beta,\mu}(\lambda;\wu,\wv)|_{\wv_L = f \wu_L},
$$
which is equivalent to the following corollary.

\begin{corollary}[single framed outer brane]
Let $(L,f)$ be a framed outer brane in $X$, let
$z$ be the unique $\bT'$-fixed point in $V(\tau_L)$, and
let $\phi \in H_{\bT'}^6(X;\bZ)$ be the $\bT'$-equivariant Poincar\'{e}
dual of $z$.   For any effective curve class  $\beta\in H_2(X;\bZ)$ and any
partition $\mu=(\mu_1,\ldots, \mu_h)$, we have
\begin{equation}\label{eqn:single-outer}
\begin{aligned}
N_{g,\beta+|\mu|b ,\mu}^{X,L}(f)=& \frac{1}{ |\Aut(\mu)| }\prod_{j=1}^h
(-1)^{f\mu_j} \frac{\prod_{m=1}^{\mu_j-1}(f\mu_j  +  m) }{\mu_j \cdot \mu_j!} \\
& \cdot \biggl( \int_{[\Mbar_{g,h}(X,\beta)^{\bT'}]^{\vir}}\frac{1}{e_{\bT'}(N^\vir)}
\prod_{j=1}^h\frac{\ev_j^*\phi}{\frac{\wu_L}{\mu_j}(\frac{\wu_L}{\mu_j}-\psi^i_j)}
\biggr)\biggr|_{\wv_L=f\wu_L}
\end{aligned}
\end{equation}
\end{corollary}

\subsection{Open GW invariants of $(X,L)$: single outer or inner brane}\label{sec:single}

Let $(L,f)$ be a framed (inner or outer) Aganagic-Vafa A-brane.
Then $L=L_1$. Let
$$
b=b_1\in H_2(X,L;\bZ), \quad \gamma=\gamma_1\in H_1(L;\bZ),\quad
x^+=x_1,
$$
where $b_1, \gamma_1,x_1$ are defined as in Section \ref{sec:topology}
Let $\phi_L^+ \in H_{\bT'}^6(X;\bZ)$ be the $\bT'$-equivariant Poincar\'{e}
dual of the $\bT'$ fixed points $z^+ \in V(\tau_L)$. When $L$ is an inner brane,
let $x^-=x_1^-$ be the other $\bT'$ fixed point in $V(\tau_L)$, and let $\phi_L^-$ be the $\bT'$-equivariant
Poincar\'{e} dual of $x^-$.  

Given $\beta' \in H_2(X,L;\bZ)$, and $\vw=(w_1,\ldots,w_h)$, where
$w_j$ are nonzero integers, we define
$$
\Mbar_{g,\beta',\vw}:= \Mbar_{g;h}(X,L\mid \beta';w_1\gamma,\ldots, w_h\gamma)
$$
where the right hand side is the moduli space of stable maps
$u:(\Si,\partial\Si = \cup_{j=1}^h R_j)\to (X,L)$, where $\Si$ is a prestable
Riemann surface of type $(g;h)$, $R_j\cong S^1$ are connected components
of the boundary $\partial\Si$ of $\Si$, $u_*[\Si]= \beta'$, 
$u_*[R_j]=w_j \gamma$. (See \cite[Section 4]{KL} for the detailed definitions.)
Then $\bT'_\bR$ acts on $\Mbar_{g,\beta',\vw}$. 
The integral
$$
\biggl( \int_{[\Mbar_{g,\beta',\vw}^{\bT'_\bR}]^\vir }\frac{1}{e_{\bT'_\bR}(N^\vir)} \biggr) \biggr|_{\wv_L=f\wu_L}.
$$
is a rational number depending on $f\in \bZ$. It is defined
up to a sign depending on choice of orientation on the virtual tangent 
bundle of $\Mbar_{g,\beta',\vw}$.  
Define
$$
\Aut(\vw)=\{\si\in S_h\mid (w_{\si(1)},\ldots, w_{\si(h)}) = (w_1,\ldots,w_h)\}.
$$

\begin{proposition} \label{prop:open-closed}
Suppose that $\beta \in H_2(X;\bZ)$, and either
\begin{itemize}
\item $L$ is an outer brane, $w_1,\ldots, w_h$ are positive integers, or
\item $L$ is an inner brane, $w_1,\ldots, w_h$ are nonzero integers.
\end{itemize}
When $L$ is an inner brane, let  $\alpha=[V(\tau_L)]\in H_2(X;\bZ)$.
$\{1,\ldots,h\}$ is a disjoint union of $J_+$ and $J_-$, where
$J_\pm=\{j\in \{1,\ldots,h\}\mid \pm w_j>0\}$. (So $J_-$ is empty
when $L$ is an outer brane.) Define $\beta' \in H_2(X,L;\bZ)$ by
$$
\beta'= \begin{cases} 
\beta + \left(\sum_{j=1}^h w_j\right)b & \textup{ if $L$ is an outer brane},\\
\beta + \left(\sum_{j\in J_+} w_j\right)b  + \left(\sum_{j\in J_-} (-w_j)\right)(\alpha-b)
& \textup{ if $L$ is an inner brane}.
\end{cases}
$$
Then
$$
\frac{1}{|\Aut(\vw)|}\int_{[\Mbar_{g,\beta',\vw}^{\bT'_\bR}]^\vir }\frac{1}{e_{\bT'_\bR}(N^\vir)} \Bigr|_{\wv_L=f\wu_L}
=\pm N_{g,\beta'\vw}(f)
$$
where
\begin{equation} \label{eqn:open-closed}
\begin{aligned}
N_{g,\beta',\vw}(f)=& \frac{1}{|\Aut(\vw)|}\prod_{j\in J_+}(-1)^{f w_j}\frac{\prod_{m=1}^{w_j-1}(fw_j + m)}{w_j\cdot w_j!}\\
& \cdot \prod_{j\in J_-}(-1)^{(f+n)w_j}\frac{\prod_{m=1}^{-w_j-1}( (f+n)(-w_j)+m)}{(-w_j)\cdot (-w_j)!} \\
&\quad \cdot\biggl( \int_{[\Mbar_{g,h}(X,\beta)^{\bT'}]^\vir} \frac{1}{e_{\bT'} (N^\vir)}\cdot 
\frac{ \prod_{j\in J_+}\ev_j^* \phi_L^+\prod_{j\in J_-}\ev_j^*\phi_L^-}{\prod_{j=1}^h[\frac{\wu}{w_j}(\frac{\wu}{w_j}-\psi_j)]}  
\biggr) \biggr|_{\wv_L=f\wu_L}
\end{aligned}
\end{equation}
\end{proposition}

\begin{remark}
The above formula \eqref{eqn:open-closed} agrees with \eqref{eqn:single-outer}
when $L$ is an outer brane. We use \eqref{eqn:open-closed} to extend the definition 
of open GW invariants when $(L,f)$ is an inner brane. 
\end{remark}

\begin{proof}[Proof of Proposition \ref{prop:open-closed}]  
When all the $w_j$'s are positive, the computation
for an inner brane is the same as that for an outer brane. From now
on, we assume that $L$ is an inner brane. $V(\tau_L)$ is the union of 
two disks $D_+$ and $D_-$ which contain the torus fixed points $x^+$ and $x^-$, respectively.
$D_+\cap D_- = L\cap V(\tau_L)$.

\begin{figure}[h]
\psfrag{u}{\small $\wu_L$} 
\psfrag{fu}{\small $f\wu_L$}
\psfrag{-(f+1)u}{\small $-(f+1)\wu_L$}
\psfrag{-u}{\small $-\wu_L$}
\psfrag{(f+n+1)u}{\small $(f+n+1)\wu_L$}
\psfrag{-(f+n)u}{\small $-(f+n)\wu_L$}
\psfrag{p}{\small $x^+$}
\psfrag{q}{\small $x^-$}
\includegraphics[scale=0.6]{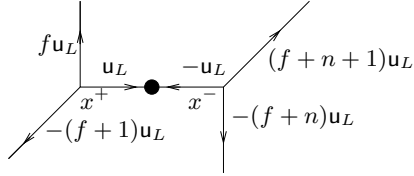}
\caption{weights of the action by $T_{L,f}$}
\label{fig:inner}
\end{figure} 

Suppose that $u:(C,q_1,\ldots, q_h)\to X$ is a stable map which represents a 
$
\bT'$ fixed point in $\Mbar_{g,h}(X,\beta)$, and that $f(q_j)=x^\pm$ for $j\in J_\pm$. Let
$\Si$ be the prestable bordered Riemann surface
$$
\Si = C\cup D_1\cup \cdots \cup D_h,\quad
$$
where $D_j$ intersects $C$ at a node $q_j$. If $j\in J_\pm$ then  
we define $u_j:D_j\to D_\pm$ by $z\mapsto z^{\pm w_j}$. Then
we obtain a stable map $\hu:(\Si,\bSi)\to (X,L)$ which represents a $\bT'_\bR$-fixed point in 
$\Mbar_{g,\beta',\vw}$. Conversely, all $\bT'_\bR$ fixed points in $\Mbar_{g,\beta',\vw}$ arise in this way.
Let  
\begin{eqnarray*}
&& \hB_1=\Aut(\Si), \  \hB_2 = \Def(\hu),\  \hB_4 = \Def(\Si),\  \hB_5=\Obs(\hu),\\
&& B_1=\Aut(C,q_1,\ldots, q_h),\  B_2 =\Def(u)\\
&& B_4 =\Def(C,q_1,\ldots,q_h),\  B_5 =\Obs(u).
\end{eqnarray*}
(See \cite[Section 4]{KL} for the definitions of these real vector spaces.) 
Let $\hT^1-\hT^2$ be the virtual tangent space to $\Mbar_{g,\beta',\vw}$ at the moduli point 
$[\hu:(\Si,\bSi)\to (X,L)]$, and let $T^1-T^2$ be the virtual tangent space to $\Mbar_{g,h}(X,\beta)$ 
at the moduli point $[u:(C,q_1,\ldots,q_h)\to X]$. Let $e=e_{\bT'_\bR}$ denote
the $\bT'_\bR$-equivariant Euler class. Then
$$
\frac{e(\hT^{2,m})}{e(\hT^{1,m})} =\frac{e(\hB^m_1) e(\hB^m_5)}{e(\hB^m_4)e(\hB^m_2)},\quad
\frac{e(T^{2,m})}{e(T^{1,m})} =\frac{e(B^m_1) e(B^m_5)}{e(B^m_4)e(B^m_2)}
$$

We have
$$
\hB_1^m = B_1^m,\quad \hB_4^m = B_4^m \oplus \bigoplus_{j=1}^h T_{q_j}C\otimes T_{q_j} D_j.
$$
Therefore
\begin{equation}\label{eqn:domain}
\frac{e(\hB_1^m)}{e(\hB_4^m)} =\frac{e(B_1^m)}{e(B_4^m)}\prod_{j\in J_+}(\frac{\wu_L}{w_j}-\psi_j)
\prod_{j\in J-}(\frac{-\wu_L}{-w_j}-\psi_j) =
\frac{e(B_1^m)}{e(B_4^m)}\prod_{j=1}^h(\frac{\wu_L}{w_j}-\psi_j)
\end{equation}

We also have a long exact sequence
\begin{eqnarray*}
0& \to \hB_2\to B_2 \oplus \bigoplus_{j=1}^h H^0(D_j)\to \bigoplus_{j\in J_+} T_{x^+} X\oplus \bigoplus_{j\in J_-} T_{x^-} X \\
 &\to \hB_5\to B_5 \oplus \bigoplus_{j=1}^h H^1(D_i)\to 0. 
\end{eqnarray*}
where
$$
H^p(D_j)= H^p(D_j,\partial D_j, u_j^*TX, (u_j|_{\partial D_j})^* TL),\quad p=0,1. 
$$
We have
\begin{eqnarray*}
e_{\bT'}(T_{x^+} X) &=& \wu_L\wv_L(-\wu_L-\wv_L)= \ev_j^*\phi^+_L \quad j\in J_+, \\
e_{\bT'}(T_{x^-} X) &=& (-\wu_L)(-\wv_L-n\wu_L)(\wv_L+(n+1)\wu) = \ev_j^*\phi^-_L \quad j\in J_-  \\
\frac{e(H^1(D_j)^m)}{e(H^0(D_j)^m)} &=& \pm \frac{1}{\wu_L}\prod_{m=1}^{w_j-1}\frac{w_j\wv_L+m\wu_L}{m\wu_L},\quad j\in J_+ \\
\frac{e(H^1(D_j)^m)}{e(H^0(D_j)^m)} &=& \pm \frac{1}{\wu_L}\prod_{m=1}^{w_j-1}\frac{w_j\wv_L+(n+m)\wu_L}{m\wu_L},\quad j\in J_-  \\
\end{eqnarray*}
Therefore
\begin{equation}\label{eqn:map}
\frac{e(\hB_5^m)}{e(\hB_2^m)}
=\frac{e(B_5^m)}{e(B_2^m)} 
\cdot\prod_{j\in J_+} \frac{\ev_j^*\phi_+^L}{\wu_L} \prod_{m=1}^{w_j-1}\frac{w_j\wv_L+m\wu_L}{m\wu_L}
\cdot\prod_{j\in J_-} \frac{\ev_j^*\phi_-^L}{-\wu_L} \prod_{m=1}^{w_j-1}\frac{w_j\wv_L+(n+m)\wu_L}{m\wu_L}
\end{equation}
Finally, 
\begin{equation}\label{eqn:aut}
\Aut(\hu)=\Aut(u)\times \prod_{j=1}^h (\bZ/w_j\bZ).
\end{equation}
Combining \eqref{eqn:domain}, \eqref{eqn:map}, and \eqref{eqn:aut}, we get
\begin{eqnarray*}
N_{g,\beta',\vw}(f)= \pm \frac{1}{|\Aut(\vw)|}\prod_{j\in J_+}\frac{\prod_{m=1}^{w_j-1}(fw_j + m)}{w_j\cdot w_j!}
\prod_{j\in J_-}\frac{\prod_{m=1}^{-w_j-1}( (f+n)(-w_j)+m)}{(-w_j)\cdot (-w_j)!} \\
\quad \cdot \biggl( \int_{[\Mbar_{g,h}(X,\beta)^{\bT'}]^\vir} \frac{1}{e_{\bT'} (N^\vir)}\cdot 
\frac{ \prod_{j\in J_+}\ev_j^* \phi_L^+\prod_{j\in J_-}\ev_j^*\phi_L^-}{\prod_{j=1}^h[\frac{\wu_L}{w_j}(\frac{\wu_L}{w_j}-\psi_j)]}  
\biggr) \biggr|_{\wv_L=f\wu_L}
\end{eqnarray*}
\end{proof}

\subsection{Generating functions for a fixed topological type of the domain} \label{sec:generating}
Any $\beta'\in H_2(X,L;\bZ)$ is of the form
$\beta'=\beta+wb$, where $\beta\in H_2(X;\bZ)$ and $w\in \bZ$. As in Section \ref{sec:geometry}, we fix a choice of $\bZ$-basis $l^{(1)},\dots, l^{(k)}$ and its dual $p_1,\dots, p_k$.
Define
$$
Q^{\beta'} = Q_0^w\prod_{a=1}^k Q_a^{d_a},
\quad d_a=\langle \beta,p_a \rangle.
$$
We define generation functions
\begin{eqnarray*}
F_{g,\vw}(Q;f) &=&\sum_{\beta'\in H_2(X,L;\bZ)} N_{g,\beta', \vw}(f) Q^{\beta'},\\
F^{g,h}(Q,Y;f) &=& \sum_{w_j\in \bZ-\{0\}} F_{g,w_1,\ldots,w_h}(Q,f)\prod_{j=1}^h Y_j^{w_j},
\end{eqnarray*}
where
$$
Q=(Q_0,Q_1,\ldots, Q_k),\quad Y=(Y_1,\ldots, Y_h).
$$

Note that when $h=1$, $w=w_1 \in \bZ$ is determined by $\partial \beta' = w\gamma$, 
so we may omit the variable $Y_1$. In particular, we define
\begin{equation}\label{eqn:FQf}
F(Q;f):= F^{0,1}(Q,Y_1=1;f) =\sum_{w\in \bZ, w\neq 0, \beta\in H_2(X;\bZ)}N_{0,\beta+wb, w}Q_0^w\prod_{a=1}^k Q_a^{d_a} ,
\end{equation}
where $d_a=\langle \beta,p_a\rangle$.

For later convenience, we introduce the following generating function for any nonzero integer $w$:
\begin{equation}\label{eqn:JL}
J^L_w(Q;f) =\begin{cases}
1+\displaystyle{ \sum_{\beta>0} Q^\beta 
\biggl(\int_{[\Mbar_{0,1}(X,\beta)^{\bT'}]^\vir} \frac{1}{e_{\bT'} (N^\vir)}
\frac{\ev_j^* \phi_L^+}{\frac{\wu_L}{w}(\frac{\wu_L}{w}-\psi_j)]}\biggr) \biggr|_{\wv_L=f\wu_L} } & w>0 \\
1+\displaystyle{ \sum_{\beta>0} Q^\beta 
\biggl( \int_{[\Mbar_{0,1}(X,\beta)^{\bT'}]^\vir} \frac{1}{e_{\bT'} (N^\vir)}
\frac{\ev_j^* \phi_L^-}{\frac{\wu_L}{w}(\frac{\wu_L}{w}-\psi_j)]}\biggr) \biggr|_{\wv_L=f\wu_L} } & w<0
\end{cases} 
\end{equation}
where $\beta\in H_2(X;\bZ)$. We say $\beta\geq 0$ when $\beta$ is effective, and $\beta>0$ when it is also non-zero. The following is an immediate consequence of Proposition \ref{prop:open-closed}.
\begin{corollary}\label{cor:FJ}
\begin{enumerate}
\item[(a)]
If $L$ is an outer brane then
\begin{equation}
F(Q;f)=\sum_{w>0}Q_0^w (-1)^{fw}\frac{\prod_{m=1}^{w-1}(fw+m)}{w\cdot w!}J_w^L(Q;f).
\end{equation}
\item[(b)] 
If is an inner brane then
\begin{equation}
\begin{aligned}
F(Q;f)=&\sum_{w>0}Q_0^w (-1)^{fw}\frac{\prod_{m=1}^{w-1}(fw+m)}{w\cdot w!}J_w^L(Q;f) \\
& +\sum_{w<0}(\frac{Q_0}{Q^\alpha})^w (-1)^{(f+n)w} \frac{\prod_{m=1}^{-w-1}((f+n)(-w)+m)}{-w\cdot(-w)!}J_w^L(Q;f)
\end{aligned}
\end{equation}
\end{enumerate}
\end{corollary}

\section{The Mirror Conjecture on Holomorphic Disks}
\label{sec:conj}

Aganagic-Vafa and Aganagic-Klemm-Vafa \cite{AV, AKV} construct a B-model Landau-Ginzburg model together with the superpotential $W$ as the mirror to the A-model on the Calabi-Yau toric threefold $X$. The superpotential $W$ is conjectured, up to a mirror transform, to be equal to the single brane disk amplitude constructed in Section \ref{sec:generating}. As discussed in Section \ref{sec:single} we fix a framed Aganagic-Vafa A-brane $(L,f)$, and let $x^+$ be the $\bT$-fixed point in $V(\tau_L)$ (when $L$ is an inner brane there is another $\bT$-fixed point $x^-$). We fix indices $i_1,i_2,i_3$, and in case of inner brane $L$, also an index $i_4$ as in Section \ref{sec:framing} (Figure \ref{fig:inner}).

\subsection{The Hori-Vafa mirror}
\label{sec:hori-vafa-mirror}
Applying $\Hom(-, \bC^*)$ to \eqref{eqn:T}, we obtain 
\begin{equation}
1\to \bT^\vee \to \tbT^\vee \stackrel{\phi}{\to} G^\vee\to 1,
\end{equation}
where $\bT^\vee \cong (\bC^*)^n$, $\tbT^\vee\cong (\bC^*)^r$,
$G^\vee \cong (\bC^*)^k$. Applying $\Hom(-,\bC^*)$ to 
\eqref{eqn:Tzero}, we obtain
\begin{equation}
1\to (\bC^*)^\vee \to \bT^\vee \to (\bT')^\vee \to 1,
\end{equation}
where  $(\bT')^\vee \cong (\bC^*)^{n-1}$.

Given $q= (q_1,\ldots,q_k)\in G^\vee$, define
$$
\bT^\vee_q :=\phi^{-1}(q)=\{\vx=(x_1,\ldots, x_{k+3})\in \tbT^\vee\mid 
\prod_{i=1}^{k+3} x_i^{l_i^{(a)}} = q_a,\ a=1,\ldots,k \} \cong (\bC^*)^3. 
$$
Let $\bC^*$ acts on $\tbT^\vee$ by $\lambda\cdot (x_1,\ldots, x_{k+3})
= (\lambda x_1,\ldots, \lambda x_{k+3})$. Then 
the $\bC^*$-action preserves $\bT^\vee_q$. Let $\bT'^\vee_q:= \bT^\vee_q/\bC^*\cong (\bC^*)^2$. 

Let $H(\vx, q)$ be the restriction of $\sum_{i=1}^{k+3} x_i$ to $\bT^\vee_q$. 
There is a $\bC^*$ action on the following Calabi-Yau $4$-fold
$$
\tX^\vee_q=\{ (w^+, w^-, \vx) \in \bC^2\times \bT^\vee_q \mid w^+ w^- = H(\vx,q)\},
$$ 
given by
\begin{eqnarray*}
\lambda\in \bC^*:  & w^-  & \mapsto \lambda w^-,\\
\ &(x_1,\dots, x_{k+3})&\mapsto (\lambda x_1,\dots,\lambda x_{k+3}).
\end{eqnarray*}
The mirror of $X$ is the quotient noncompact Calabi-Yau 3-fold $\tX^\vee_q/\bC^*$, denoted by $X^\vee_q$. 
The curve in $X^\vee_q$ given by the equation $H(\vx,q)=0$, $w^+=w^-=0$ is the \emph{mirror curve}.

For $a=1,\ldots, k$, we extend the charge vectors $l^{(a)}=(l^{(a)}_1,\ldots, l^{(a)}_{k+3})$ to
$$
\tl^{(a)} = (l^{(a)}_1,\ldots, l^{(a)}_{k+3} , 0, 0).
$$
There is an additional charge vector $\tl^{(0)}=(l^{(0)}_1,\ldots, l^{(0)}_{k+5})$ for $i=1,\dots, k+5$ with 
$$
\sum_{i=1}^{k+5}l^{(0)}_i=0
$$ that characterizes the B-brane. Define the index sets
\begin{equation}
\label{eqn:I}
I=\{1,\dots, r\},\quad  \tI=\{1,\dots, r+2\},\quad
I_0=\{i_1,i_2,i_3\}.
\end{equation}

If the A-brane carries the framing $f$, we set the open charge vector $l^{(0)}=(l^{(0)}_1,\ldots, l^{(0)}_{k+3})$ 
to be 
$$
l^{(0)}_{i_1}= 1,\quad l^{(0)}_{i_2}=f, \quad
l^{(0)}_{i_3}= -f-1,\quad l^{(0)}_i =0 \quad \textup{for }i\in I\backslash I_0.
$$
We introduce an extra charge vector 
$$
\tl^{(0)}=(l^{(0)}_1,\ldots, l^{(0)}_r,1,-1).
$$
The mirror curve can be written down in the following coordinate patch 
$$x_{i_1}=-\tx,\ x_{i_2}=\ty,\ x_{i_3}=1.$$
Other $x_i$ for $i\in I$ are determined from $\tx,\ty$ through the equation $$\prod_{i=1}^{k+3} x_i^{l_i^{(a)}} =q_a,\quad a=1,\ldots,k .$$The mirror curve is an affine curve in $\bC^*\times \bC^*$ given by a single equation $$H(\tx,q)=-\tx+\ty+1+\sum_{i\in I\backslash I_0} x_i(\tx,\ty, q)=0$$ in coordinates $\tx, \ty$.

\subsection{Aganagic-Vafa B-branes}
As shown in \cite{AV}, the mirror B-brane of the A-brane defined in Section \ref{sec:A-brane} above is 2-cycle in $X^\vee_q$, 
given by the equations
\begin{equation}
\label{eqn:b-brane}
w^-=0,\ H(\vx,q)=0,
\end{equation}
where $\vx=\vx(w^+)\in \{\vx|H(\vx,q)=0\}$ is a function of the coordinate $w^+$ on the brane. We fix $\vx$ at $w^+=\infty$ to be $\vx(\infty)=\vx^*$. Following \cite{LM, LMW}, the framed modulus $x$ of the B-brane is determined by the value of $\vx$ on $w^{+}=0$ 
$$
x=q_0=-\prod_{i=1}^r x_i(0)^{l^{(0)}_i}=\tx \ty^f.
$$
The the B-brane modulus in framing $0$ is given by $\tx$. We denote the corresponding B-brane by $\cC_x$, and the reference B-brane by $\cC_{x^*}$, with $x^*=-\prod_{i=1}^{k+3} x_i^*(0)^{l^{(0)}_i}$. The reference brane $\cC_{x^*}$ could be chosen as a holomorphic curve with constant $x(w^+)=x^*$. Other branes cannot be holomorphic, since $x$ has different bounded values at $w^+=0$ and $w^+=\infty$.

The mirror curve in coordinates $\tx, \ty$ can be further written as an equation in the framed coordinates $x, y$ 
\begin{equation}
\label{eqn:mirror-curve}
M(x,y,q)=H(\vx,q)=-xy^{-f}+y+1+\sum_{i\in I\backslash I_0} x_i(xy^{-f},y,q)=0
\end{equation}
by a change of coordinates $$x=\tx \ty^f,\ y=\ty.$$

\subsection{Picard-Fuchs equations and the mirror map} \label{sec:PF}

Let $\Gamma(x)$ be a $3$-chain in $X^\vee_q$ such that $\partial \Gamma(x)=\cC_x-\cC_{x^*}$, and 
$\gamma\in H^3(X_q^\vee;\bZ)$ is any $3$-cycle. Mayr and Lerche-Mayr \cite{Ma, LM} show that the periods 
$$
\int_{\Ga(x)} \Omega
$$ are eliminated by  GKZ-type operators
$$
\cD_a=\prod_{l^{(a)}_i>0}\prod_{j=0}^{l^{(a)}_i-1}(\sum_{b} l^{(b)}_i\theta_b-j)
-q_a\prod_{l^{(a)}_i<0}\prod_{j=0}^{l^{(a)}_i-1}(\sum_{b=1}^k l^{(b)}_i\theta_b-j),
$$
where $a=0,\dots, k$, $i=1,\ldots,r+2$,  and $\theta_a=q^a\frac{d}{dq^a}$. The extended Picard-Fuchs equations are
\begin{align}
\label{eqn:PF}
\cD_\alpha W(q_0,q_1,\dots,q_k) = 0.
\end{align}

Let $q_a=e^{t_a}$, and the string flat coordinates $Q_a=e^{T_a}$. 
The open string variable is $Q_0=e^{T_0}$. Recall that for  $a=1,\ldots,k$, $T_a =-r_a+\sqrt{-1} \theta_a$
are complexified K\"{a}hler parameters.
%
The mirror map is
$$
T_a=\Pi^1_a(q).
$$
where $\Pi^1_a(q)$ is a solution to the Picard-Fuchs equations above, i.e. $\cD_b \Pi^1_a(q)=0$ for $b=0,\dots, k$. 
It has the leading term behavior
\begin{equation} \label{eqn:mirror-pf}
T_a= \log(q_a)+S_a(q),
\end{equation}
where $S_a(q)$ is a power series in $q_1,\dots, q_k$.

Recall that $d_a=\langle p_a, \beta \rangle$ for $\beta\in H_2(X;\bZ)$, and $\langle D_i^*, \beta \rangle =\sum_{a=1}^k d_a l^{(a)}_i$. Denote $q^\beta=\prod_{a=1}^k q_a^{d_a}$.
For any $\beta \in \bL_\eff$ and $i_0\in I$, if $\langle D_{i_0}^*,\beta \rangle <0$ and $\langle D_{i_0}^*,\beta \rangle \geq 0$ for $i\in I\backslash \{i_0\}$, we define 
\begin{align*}
E_{i_0}(\beta)= 
\frac{(-1)^{(-\langle D_{i_0}^*,\beta \rangle -1)}(-\langle D_{i_0}^*,\beta \rangle -1)!}{\prod_{i\in I\backslash \{i_0\}} \langle D_{i}^*,\beta \rangle !};
\end{align*}
otherwise we define $E_{i_0}(\beta)=0$. Recall that $\beta\geq 0$ when $\beta$ is effective, and $\beta>0$ when it is also non-zero. In particular, using the Frobenius method, as shown in \cite{M,LM}, the mirror correction is given by 
\begin{eqnarray}
S_a(q) = \sum_{i \in I} l^{(a)}_{i} A_{i}(q),\quad a=0,\ldots, k,
\label{eqn:mirror_map}
\end{eqnarray}
where $A(q)$ is the following series of $q_1,\dots, q_k$:
\begin{equation}\label{eqn:A}
A_{i}(q) =\sum_{\beta> 0} E_{i}(\beta)q^\beta.
\end{equation}
The above formula is a direct reformulation of the equations \cite[(2.20,2.21,2.22)]{LM}. 

\begin{remark}
In \cite{M,LM}, Lerche-Mayr explicitly compute the open mirror correction $S_0(q)$. Let $I'\subset \{1,\dots, r\}$ be the index set such that $A_i(q)=0$ for $i\notin I'$. One chooses an index subset $K\subset \{1,\dots,k\}$ such that $|K|=|I'|$ and $S_a$ are linearly independent. Denote the matrix $L=(l^{(a)}_{i})^{a\in K}_{i\in I'}$. Lerche-Mayr in \cite{LM} gives the open mirror correction as $S_0(q)=\sum_{i\in I'} l_{i}^{(0)} (\sum_{a\in K} (L^{-1})^{i}_a S_a(q))$. It turns out that
\begin{eqnarray*}
S_0(q)&=&\sum_{i\in I'} l_{i}^{(0)} (\sum_{a\in K} (L^{-1})^{i}_a S_a(q))\\
&=&\sum_{i\in I'} l_{i}^{(0)} (\sum_{a\in K} (L^{-1})^{i}_a \sum_{j\in I'} l^{(a)}_{j} A_{j}(q))\\
&=&\sum_{i\in I'} l_{i}^{(0)} A_{i}(q).
\end{eqnarray*}
It agrees with \eqref{eqn:mirror_map}.
\end{remark}


\subsection{Superpotential and the mirror prediction of disk invariants}

In \cite{AV}, the superpotential $W(\cC_x)$ associated to the B-brane $\cC_x$ is the integral over the chain
$$
W=\int_{\Ga(x)}\Omega,
$$
where the boundary of the chain $\partial \Ga(x)=\cC_x-\cC_{x^*}.$ 
Let $\tbL= \bZ \times \bL\cong H_2(X,L;\bZ)$. For any $\tbeta=(w,\beta) \in \tbL$, denote $\tq^\tbeta=q_0^w q^\beta$. We define the extended pairing
$$\langle D_i^*, \tbeta \rangle = w \tl_i^{(0)} + \langle D_i^*,\beta \rangle =  \sum_{a=0}^k d_a \tl_i^{(a)}.$$
Define the cone $\tbL_\eff(L)$ in $\bZ\times \bL_\eff \subset \tbL$
$$
\tbL_\eff(L)=
\{\tbeta=(w,\beta) \in \tbL|w\neq 0, \langle D_i^*,\tbeta \rangle \geq 0, i\in I\backslash\{i_2, i_3\}\}.
$$

Given the charge vectors $\tl^{(a)}$, $a=0,\ldots,k$,  Lerche and Mayr \cite{M,LM} show that $W$ is a double logarithm solution of the Picard-Fuchs equations. Precisely $W$ consists of a double logarithm part and a power series part
$$
W=\sum_{i,j=0}^k c_{ij} \log q_i \log q_j +W_0,
$$
where $W_0$ is explicitly in the following form
\begin{equation} \label{eqn:superpotential} 
W_0=
\sum_{\tbeta\in \bL_\eff(L)} C(\tbeta)  \tq^\tbeta.
\end{equation}
The coefficient $C(\tbeta)$ is obtained by applying the Frobenius method
\begin{align*}
C(\tbeta)&=
\begin{cases}
\displaystyle{
\frac{\prod_{i=i_3, r+3} (-1)^{ -\langle D_i^*, \tbeta\rangle -1}(-\langle D_{i}^*, \tbeta\rangle -1)!}
{\prod_{i\in \tI\backslash \{i_3,r+3\}}\langle D_i^*, \tbeta \rangle !} },
&  \langle D_{i_3}^*,\tbeta \rangle <0 \\
& \\
\displaystyle{
\frac{\prod_{i=i_2, r+3} (-1)^{ -\langle D_i^*, \tbeta\rangle -1}(-\langle D_{i}^*, \tbeta\rangle -1)!}
{\prod_{i\in \tI\backslash \{i_2,r+3\}}\langle D_i^*, \tbeta \rangle !} }, 
& \langle D_{i_2}^*, \tbeta \rangle <0
\end{cases}\\
&=\frac{(-1)^{\langle D_{i_3}^*+D_{r+3}, \tbeta \rangle }}{\prod_{i\in \tI\backslash \{i_2,i_3,r+3\}}\langle D_i^*, \tbeta \rangle !}\frac{\prod_{m=-\infty}^{-\langle D_{i_3}^*,\tbeta\rangle -1} m }{\prod_{m=-\infty}^{\langle D_{i_2}^*,\tbeta \rangle} m }.
\end{align*}

In \cite{AV} and \cite{AKV}, the superpotential $W_0$ is explicitly computed by solving the mirror curve $M(x,y,q)=0$. Writing $y=y(x,q)$, and in \cite{AV} it is shown that the superpotential is $$x\frac{\partial W_0}{\partial x} = -\log y(x,q)\quad  \text{up to a function in $q$}.$$
\begin{remark}
The above formula differs by a sign with Equation 4.5 of \cite{AV}. We choose this sign convention in order to match the sign convention on Gromov-Witten theory. One may simply apply a change of coordinate $y \mapsto 1/y$ to get rid of this minus sign. We only consider invariants of non-zero winding numbers. Bouchard and Su\l kowski discuss the constant contributions of the higher genus superpotential associated to a mirror curve \cite{BS}.
\end{remark}

\begin{conjecture} \label{AVconj}
After a change of variables by the mirror transform \eqref{eqn:mirror-pf},
$$
F(Q; f) = W_0(q;f).
$$
\end{conjecture}
We prove this conjecture in the next section.

\section{Proof of Conjecture  \ref{AVconj}}
\label{sec:proof}

\subsection{Equivariant mirror theorem}
A toric manifold $X$ is semi-Fano if $c_1(T_X)\geq 0$.
In this subsection, we state an equivariant mirror theorem
for a semi-projective semi-Fano toric manifold $X$ (Theorem \ref{thm:MP}). 
When $X$ is projective, or is the total space of a direct sum
of negative line bundles over a projective toric manifold $Z$,
Theorem \ref{thm:MP} follows from the results in \cite{LLY1, Gi1} (when $X$ or $Z$ is a projective space) 
and in \cite{LLY2, Gi2} (in the general case). When $X$ is a general noncompact, semi-projective,
semi-Fano toric manifold, Theorem \ref{thm:MP} follows from the results in \cite{CCIT}.

Let $X$ be a semi-projective toric manifold. 
Let $D_1,\ldots,D_r$ be the $\bT$-divisor of $X$. Let 
$$
D_i^*=c_1(\cO_X(D_i))\in H^2(X;\bZ), \quad (D_i^*)_\bT = (c_1)_\bT(\cO_X(D_i)) \in H^*_{\bT}(X;\bZ).
$$

Define an $H^*_\bT(X;\bZ)$-valued function $I(t,\hbar^{-1})$: 
$$
I(t,\hbar^{-1}) =
e^{(t_0' + \sum_{a=1}^k p_a t_a)/\hbar}\Bigl(1+ \sum_{\beta>0 } q^\beta 
\prod_{i=1}^{r}\frac{\prod_{m=-\infty}^0 ( (D^*_i)_\bT + m\hbar)} 
{\prod_{m=-\infty}^{\langle \beta, D^*_i \rangle  }( (D^*_i)_\bT + m\hbar)}\Bigr).
$$
Note that the parameter $t_0'$ is different from the parameter
$t_0=\log q_0$ introduced in Section \ref{sec:PF}.

Define another $H^*_\bT(X;\bZ)$-valued function $J(T,\hbar^{-1})$ as follows:
\begin{eqnarray*}
&& \langle \gamma , J\rangle\\
&=& \int_X e^{(T'_0+ \sum_{a=1}^k p_a T_a)/\hbar} \gamma 
 +\sum_{\beta>0} Q^\beta \int_{[\Mbar_{0,2}(X,d)]^\vir} 
\frac{\ev_1^*(e^{(T'_0 + p_1 T_1+\cdots + p_k T_k)/\hbar} \gamma)\ev_2^*(1)}{\hbar-\psi_1}\\
&=& \int_X e^{(T'_0+ \sum_{a=1}^k p_a T_a)/\hbar} \gamma
 +\sum_{\beta > 0} Q^\beta \int_{[\Mbar_{0,1}(X,d)]^\vir} 
\frac{\ev_1^*(e^{(T'_0 + \sum_{a=1}^k p_a T_a)/\hbar} \gamma)}{\hbar(\hbar-\psi_1)};
\end{eqnarray*}
or equivalently,
$$
J= e^{(T'_0+ \sum_{a=1}^k T_a p_a)/\hbar}\Bigl(1+\sum_{\beta > 0} Q^\beta (\ev_1^d)_* 
\frac{1}{\hbar(\hbar-\psi_1)}\Bigr).
$$

Assume that $X$ is semi-Fano. Then the expansion of the $I$ function in $\hbar^{-1}$ is of the following form
$$
I(t,\hbar^{-1})=1+\hbar^{-1}(t'_0 + f_0(q)+\sum_{i=1}^r \lambda_i g_i(q)  +\sum_{a=1}^k p_a(t_a+f_a(q)))+O(\hbar^{-2}),
$$
where $\lambda_1,\ldots, \lambda_r \in H^2(B\bT)$ are the universal first Chern classes.
The expansion of the $J$ function in $\hbar^{-1}$ is
$$
J(T,\hbar^{-1})=1+\hbar^{-1}(T'_0+\sum_{i=1}^k p_a T_a)+O(\hbar^{-2}).
$$
The mirror map is  given by 
$$
T_0'= t_0'+ f_0(q)+\sum_{i=1}^r \lambda_i g_i(q), \quad T_a = t_a + f_a(q), \quad a=1,\ldots,k.
$$
In particular, 
$$
Q_a=q_a\exp(f_a(q)), \quad a=1,\ldots, k.
$$ 

Theorem \ref{thm:MP} below is a special case of the results in \cite{CCIT}. 
When $X$ is the total space of a direct sum of negative line bundles
over a projective toric manifold (such as the six toric Calabi-Yau 3-folds
considered in Section \ref{sec:conifold}--\ref{sec:dP3}), 
Theorem \ref{thm:MP} is a special case of the results in \cite{LLY2, Gi2}.
(It is stated as Corollary 4.3 in \cite{Gi2}.)

\begin{theorem}[Equivariant mirror theorem]\label{thm:MP}
With the above identification 
$$
T_0'= t_0'+ f_0(q)+\sum_{i=1}^r \lambda_i g_i(q), \quad T_a = t_a + f_a(q), \quad a=1,\ldots,k,
$$ the $I$ and $J$ functions are equal 
$$
I(t,\hbar^{-1})=J(T,\hbar^{-1}).
$$
\end{theorem}

\subsection{Open and closed mirror maps}
We now restrict to the subtorus $\bT_{L,f}$. 
Let $\iota_+: x^+\to X$ be the inclusion map, and let
$\phi^+\in H^6_{\bT_{L,f}}(X;\bZ)$ be the $\bT_{L,f}$-equivariant Poincar\'{e} dual of $x^+$. 
Recall that $J_w^L(Q;f)$ is  defined by \eqref{eqn:JL} and $A_{i}(q)$ is defined by in \eqref{eqn:A}.
\begin{proposition}\label{prop:Jw}
Let $w$ be a positive integer. Then
$$
J_w^L(Q;f)
= (\frac{q_0}{Q_0})^w \left(1+ \sum_{\beta> 0} q^\beta  
\prod_{i=1}^r\frac{\prod_{m=-\infty}^0 ( w l_i^{(0)} + m)}
{\prod_{m=-\infty}^{\langle \beta, D_i^*\rangle }(w l_i^{(0)}  + m)}\right).
$$
where $Q=(Q_0,Q_1,\ldots,Q_k)$ and $q=(q_0,q_1,\ldots,q_k)$ are 
related by the following open and closed mirror maps
\begin{equation}\label{eqn:open-closed-mirror}
Q_a=q_a\exp(\sum_{i\in I} l_{i}^{(a)} A_{i}(q)),\quad  a=0,1,\ldots,k.
\end{equation}
\end{proposition}
\begin{proof} In this proof, $\wu=\wu_L$. The parameters $t_0'$ and $t_0$ are different.
$$
\iota_+^* J(T,\hbar^{-1}) = e^{(T_0' +\sum_{a=1}^k T_a \iota_+^*p_a)/\hbar}
\Bigl(1+\sum_{\beta>0} Q^\beta \int_{[\Mbar_{0,1}(X,\beta)]^\vir}\frac{\ev_1^*\phi^+}{\hbar(\hbar-\psi_1)}\Bigr).
$$
\begin{eqnarray*}
\iota_+^* J(T,\frac{w}{\wu}) &=& e^{(T'_0 +\sum_{a=1}^k T_a \iota_+^*p_a)\frac{w}{\wu}}
\Bigl(1+\sum_{\beta>0} Q^\beta \int_{[\Mbar_{0,1}(X,\beta)]^\vir}\frac{\ev_1^*\phi^+}{\frac{\wu}{w}(\frac{\wu}{w}-\psi_1)}\Bigr)\\
&=&  e^{(T'_0 +\sum_{a=1}^k T_a \iota_+^*p_a)\frac{w}{\wu}} J_w^L(Q;f)\\
\iota_+^* I(t,\frac{w}{\wu}) &=& e^{(t'_0 + \sum_{a=1}^k t_a \iota_+^* p_a)\frac{w}{\wu}} 
\Bigl(1+ \sum_{\beta>0} q^\beta  \prod_{i=1}^r\frac{\prod_{m=-\infty}^0 (l_i^{(0)}\wu + m\frac{\wu}{w})}
{\prod_{m=-\infty}^{\langle \beta,D_i^*\rangle}(l_i^{(0)}\wu + m\frac{\wu}{w})}\Bigr) \\
 &=& e^{(t_0' + \sum_{a=1}^k t_a \iota_+^* p_a)\frac{w}{\wu}} \Bigl(1+ \sum_{\beta>0} q^\beta  
\prod_{i=1}^r\frac{\prod_{m=-\infty}^0 ( w l_i^{(0)} + m)}
{\prod_{m=-\infty}^{\langle \beta,D_i^*\rangle}(w l_i^{(0)}+ m)}\Bigr) 
\end{eqnarray*}
where we used $\sum_{i=1}^r D_i^* =0$. We conclude that
$$
J^L_w(Q;f) = (\frac{q_0}{Q_0})^w \left(1+ \sum_{\beta> 0} q^\beta  
\prod_{i=1}^r\frac{\prod_{m=-\infty}^0 ( w l_i^{(0)} + m)}
{\prod_{m=-\infty}^{\langle \beta, D_i^*\rangle }(w l_i^{(0)}  + m)}\right).
$$
where $Q_0$ and $q_0$ are related by the following open mirror map
$$
Q_0=q_0\exp( \frac{T_0'-t_0'}{\wu} + \sum_{a=1}^k (T_a-t_a)\frac{\iota_+^* p_a}{\wu}), 
$$
or equivalently,
\begin{equation}\label{eqn:Ttzero}
T_0 = t_0 + \frac{ \iota_+^* ( (T_0'-t_0') + \sum_{a=1}^k (T_a-t_a)p_a) }{\wu}.
\end{equation}

It remains to show that the open and closed mirror maps are given by 
$$
T_a = t_a + \sum_{i\in I} {l^{(a)}_i} A_i(q), \quad a=0,1,\ldots,k.
$$
To see this, we expand the $I$ function from Section \ref{sec:mirror-principle}:
\begin{eqnarray*}
I(t,\hbar^{-1})&=&(1+\frac{t_0'+\sum_{a=1}^k t_ap_a}{\hbar}+O(\hbar^{-2}))(1+\sum_{\beta > 0} q^\beta 
\prod_{i=1}^r \frac{\prod_{m=-\infty}^0( (D^*_i)_\bT +m\hbar)}{\prod_{m=-\infty}^{\langle\beta,D_i \rangle}((D_i^*)_\bT+m\hbar)})\\
&=&(1+\frac{t_0'+\sum_{a=1}^k t_a p_a }{\hbar}+O(\hbar^{-2})) (1+ \sum_{i\in I} \hbar^{-1} A_{i}(q)(D_{i}^*)_\bT +O(\hbar^{-2})).
\end{eqnarray*}
The $\hbar^{-1}$-term in $I(t,\hbar^{-1})$ is
$$
\hbar^{-1}\text{-term}=\hbar^{-1}((t'_0+\sum_{a=1}^k t_ap_a)+ \sum_{i\in I} A_{i}(q) (D_{i}^*)_\bT)
$$
Therefore, the mirror maps satisfy
\begin{equation}\label{eqn:equivariant}
(T_0'-t_0')+ \sum_{a=1}^k (T_a-t_a) p_a = \sum_{i\in I} A_{i} (q) (D_{i}^*)_\bT
\end{equation}

To obtain the closed mirror map, we consider the nonequivariant
version of \eqref{eqn:equivariant}:
$$
f_0(q)+ \sum_{a=1}^k (T_a-t_a) p_a = \sum_{i\in I} A_{i}(q) D_{i}^*
$$
where $D_{i}^*= \sum_{a=1}^k l^{(a)}_{i} p_a \in H_2(X;\bZ)$. So $f_0(q)=0$, and
the closed mirror map is given by
$$
T_a = t_a +\sum_{i\in I} l^{(a)}_{i} A_{i}(q),\quad a=1,\ldots,k.
$$

To obtain the open mirror map, we apply $\iota_+^*$ to \eqref{eqn:equivariant} and restrict to $\bT_{L,f}$:
\begin{equation}\label{eqn:iota-p}
\iota_+^* \left( (T_0'-t_0') + \sum_{a=1}^k (T_a-t_a)p_a \right) = \sum_{i\in I} A_{i}(q) \iota_+^* (D_{i}^*)_{\bT_{L,f}} = \sum_{i\in I} A_{i}(q)l_{i}^{(0)}\wu.
\end{equation}
Combing \eqref{eqn:Ttzero} and \eqref{eqn:iota-p}, we obtain the open mirror map:
$$
T_0=t_0 + \sum_{i\in I} l^{(0)}_{i} A_{i}(q).
$$
\end{proof}

Note that \eqref{eqn:open-closed-mirror} implies
\begin{equation}
Q^\beta =q^\beta \exp\left(\sum_{i\in I} \langle \beta,D_{i}^*\rangle A_i(q)\right).
\end{equation}
for all $\beta\in H_2(X;\bZ)$.

\begin{remark}
The mirror maps \eqref{eqn:open-closed-mirror} are exactly prescribed by the logarithmic solution of the Picard-Fuchs equations \eqref{eqn:mirror-pf}, as expected.
\end{remark}

\begin{corollary}
\label{cor:Jwminus}
Let $w$ be a negative integer, and $L$ be a framed inner brane. Define $l^- = (l^-_1,\ldots, l^-_{k+3})$
to be
$$
l^-_{i_2}= f+n+1,\quad l^-_{i_3}=-(f+n), \quad l^-_{i_4}= -1,
$$
and
$$
l^- =0 \textup{ if } i\in \{1,2,\ldots,k+3\}-\{i_2, i_3, i_4\}.
$$
(In other words, $l^-_i =l^{(0)}_i -\langle D_i^*, \alpha \rangle,\quad i=1,\ldots, k+3.$) Then  $$J^L_w(Q;f)=(\frac{q_0/q^\alpha}{Q_0/Q^\alpha})^w
\left(1+ \sum_{\beta>0} q^\beta \prod_{i=1}^{k+3}\frac{\prod_{m=-\infty}^0 (w l^-_i + m)}
{\prod_{m=-\infty}^{\langle \beta,D_i^*\rangle } (w l^-_i + m)} \right).$$
\end{corollary}
\begin{proof}
 
Apply Proposition \ref{prop:Jw} with respect to the fixed point $q$ and $L$, we have
$$J^L_w(Q;f)=(\frac{q_0^-}{Q_0^-})^w
\left(1+ \sum_{\beta>0} q^\beta \prod_{i=1}^{k+3}\frac{\prod_{m=-\infty}^0 (w l^-_i + m)}
{\prod_{m=-\infty}^{\langle \beta,D_i^*\rangle } (w l^-_i + m)} \right),$$
where $q_0^-$ and $Q_0^-$ are open string parameters related by 
$$
Q_0^-=q_0^-\exp(\sum_{i\in I} l^-_i A_i(q)).
$$ 
Since $l^-_i =l^{(0)}_i -\langle \alpha, D_i^*\rangle$, we have
$$
Q_0^-=q_0^-\exp\Big(\sum_{i\in I}l^{(0)}_i A_i(q)-\sum_{i\in I} \langle \alpha, D_i^* \rangle A_i (q) \Big).
$$
On the other hand,
$$
\frac{Q_0/Q^\alpha}{q_0/q^\alpha}=
\exp\Big(\sum_{i\in I}l^{(0)}_i A_i(q)-\sum_{i\in I} \langle \alpha, D_i^* \rangle A_i(q) \Big)
$$
and this implies the statement.
\end{proof}

\subsection{Outer brane}
By Corollary \ref{cor:FJ} (a) and Proposition \ref{prop:Jw}, 
\begin{align*}
F(Q;f)
&= \sum_{w>0}q_0^w(-1)^{f w} \frac{\prod_{m=1}^{w-1}(f w+m)}{w \cdot w!}
\left(1+ \sum_{\beta>0} q^\beta \prod_{i=1}^{k+3}\frac{\prod_{m=-\infty}^0 (w l^{(0)}_i + m)}
{\prod_{m=-\infty}^{\langle \beta,D_i^*\rangle } (w l^{(0)}_i + m)} \right) \\
&= \sum_{(w,\beta)\in\tbL_\eff(L), w>0}
(-1)^{f w} \frac{\prod_{m=1}^{w-1}(f w+m)}{w \cdot w!} \cdot \frac{1}{\prod_{i\in I\backslash I_0} \langle D_i^*,\beta \rangle!} \\
& \cdot \frac{\prod_{m=-\infty}^{0} (w+m) }{\prod_{m=-\infty}^{\langle D_{i_1}^*,\beta \rangle}(w+m) } \cdot \frac{\prod_{m=-\infty}^0 (fw + m)} {\prod_{m=-\infty}^{\langle D_{i_2}^*, \beta \rangle } (fw + m)}
\cdot \frac{\prod_{m=-\infty}^0 (-(f+1)w + m)} {\prod_{m=-\infty}^{\langle D_{i_3}^*,\beta \rangle } (-(f+1)w + m)} q_0^w q^\beta\\
&=\sum_{\tbeta\in \tbL_\eff(L)}  \frac{(-1)^{\langle D_{i_3}^*+D_{r+3}^*,\tbeta\rangle}}{\prod_{i\in \tI\backslash \{i_2,i_3\}} \langle D_i^*,\tbeta \rangle!} \frac{\prod_{m=-\infty}^{ -\langle D_{i_3}^*,\tbeta \rangle -1} m }{\prod_{m=-\infty}^{ \langle D_{i_2}^*, \tbeta \rangle} m }\tq^\tbeta.
\end{align*}
This agrees with Equation \eqref{eqn:superpotential}.
\begin{remark}
The condition $\tbeta=(w,\beta) \in \tbL_\eff(L)$ implies
$$
w\neq 0, \quad \langle D_{i_1}^*,\tbeta \rangle = w+ \sum_{a=1}^k d_a l^{(a)}_{i_1} \geq 0.
$$
When $L$ is an outer brane, it further implies $w>0$. One may choose appropriate $\bZ$-basis $\{l^{(a)}\}_{a=1,\dots,k}$ in $\bL$ and $\{p_a\}_{a=1,\dots,k}$ in $\bL^\vee$ with respect to the $\bT$-fixed point $x^+$. Since the toric variety $X=(\bC^r-Z(\Delta))/G$ is smooth, the $G$ action on $\bC^r-Z(\Delta)$ is free; in particular the $G$ action on 
$$\{(X_1,\dots, X_r): X_{i_1}=X_{i_2}=X_{i_3}=0, X_i \in \bC^* \text{ for other } i\}\subset \bC^r-Z(\Delta)$$ 
is free. It follows that the composition with the projection
$$
\bL \to \tN=\bigoplus_{i=1}^r \bZ \tv_i \to \bigoplus_{i\in I\backslash I_0} \bZ \tv_i
$$
is an isomorphism. Thus we choose the preimages of $\tv_i$, where $i\in I\backslash I_0$ as a basis $\{l^{(a)}\}$ of $\bL$. We write $l^{(a_i)}$ to be the preimage of $\tv_i$ for $i \in I\backslash I_0$. It is obvious that the cone $$\sum_{a=1}^k \bZ_{\geq 0} l^{(a)} \subset \bL_\eff;$$
while the pairing $\langle p_a, l^{(b)} \rangle= \delta_{ab}$ for the dual basis $p_a\in \bL^\vee$, and $p_{a_i}=D_i^*$. With such a choice of $l^{(a)}$ and $p_a$, the intersection $L\cap X$ is described by the following equations
\begin{align*}
l^{(a)}_{i_1} |X_{i_1}|^2+\sum_{i\in I\backslash I_0} l^{(a)}_i |X_i|^2&=r_a,\quad a=1,\dots,k,
\end{align*}
where $r_1,\dots,r_k$ are K\"ahler parameters. Since $L$ is an outer brane, if one moves $L$ to infinity, i.e. let $|X_{i_1}|\to \infty$, we will not encounter a $\bT$-fixed point, i.e. $|X_i|$ cannot be zero for any $i\in I\backslash I_0$. By our choice of basis $\{l^{(a)}\}$ and $\{p_a\}$ in $\bL$ and $\bL^\vee$, $l^{(a_i)}_j=\delta_{ij}$ for $i,j\in I\backslash I_0$. Thus $l_{i_1}^{(a)}\leq 0$ for $a=1,\dots, k$, and $\tbeta\in\tbL_\eff(L)$ already implies $w>0$.
\end{remark}

This completes the proof of Conjecture \ref{AVconj} when $L$ is a framed outer brane.

\subsection{Inner brane}
By Corollary \ref{cor:FJ} (b), Proposition \ref{prop:Jw} and Corollary \ref{cor:Jwminus},
$$
F(Q;f) = I_+ + I_-,
$$
where
\begin{eqnarray*}
I_+&=& \sum_{w>0}q_0^w(-1)^{f w} \frac{\prod_{m=1}^{w-1}(f w+m)}{w \cdot w!}
\left(1+ \sum_{\beta>0} q^\beta \prod_{i=1}^{k+3}\frac{\prod_{m=-\infty}^0 (w l^{(0)}_i + m)}
{\prod_{m=-\infty}^{\langle \beta,D_i^*\rangle } (w l^{(0)}_i + m)} \right) \\
I_-&=& \sum_{w<0}(\frac{q_0}{q^\alpha})^w(-1)^{(f+n) w} \frac{\prod_{m=1}^{-w-1}((f+n)(-w)+m)}{(-w) \cdot(-w)!}
\left(1+ \sum_{\beta>0} q^\beta \prod_{i=1}^{k+3}\frac{\prod_{m=-\infty}^0 (w l^-_i + m)}
{\prod_{m=-\infty}^{\langle \beta,D_i^*\rangle } (w l^-_i + m)} \right).
\end{eqnarray*}
So $I_-$ can be rewritten as 
\begin{eqnarray*}
I_- &=& \sum_{w<0}q_0^w(-1)^{(f+n) w} \frac{\prod_{m=1}^{-w-1}((f+n)(-w)+m)}{(-w) \cdot(-w)!}
\left(\sum_{\beta\geq 0} q^{\beta-w\alpha} \prod_{i=1}^{k+3}\frac{\prod_{m=-\infty}^{-w\langle \alpha,D_i^*\rangle} (w l^{(0)}_i + m)}
{\prod_{m=-\infty}^{\langle \beta-w\alpha,D_i^*\rangle } (w l^{(0)}_i + m)} \right) \\
&=& \sum_{w<0}q_0^w(-1)^{(f+n) w} \frac{\prod_{m=1}^{-w-1}((f+n)(-w)+m)}{(-w) \cdot(-w)!}
\left(\sum_{\beta+w\alpha\geq 0} q^\beta \prod_{i=1}^{k+3}\frac{\prod_{m=-\infty}^{-w\langle \alpha,D_i^*\rangle} (w l^{(0)}_i + m)}
{\prod_{m=-\infty}^{\langle \beta,D_i^*\rangle } (w l^{(0)}_i + m)} \right) \\\end{eqnarray*}
The calculation of $I_+$ is identical to the outer brane case.
\begin{align*}
I_+
&= \sum_{w>0}q_0^w(-1)^{f w} \frac{\prod_{m=1}^{w-1}(f w+m)}{w \cdot w!}
\left(1+ \sum_{\beta>0} q^\beta \prod_{i=1}^{k+3}\frac{\prod_{m=-\infty}^0 (w l^{(0)}_i + m)}
{\prod_{m=-\infty}^{\langle \beta,D_i^*\rangle } (w l^{(0)}_i + m)} \right) \\
&=\sum_{\tbeta\in \tbL_\eff(L), w>0}  \frac{(-1)^{\langle D_{i_3}^*+D_{r+3}^*,\tbeta\rangle}}{\prod_{i\in \tI\backslash \{i_2,i_3\}} \langle D_i^*,\tbeta \rangle!} \frac{\prod_{m=-\infty}^{ -\langle D_{i_3}^*,\tbeta \rangle -1} m }{\prod_{m=-\infty}^{ \langle D_{i_2}^*, \tbeta \rangle} m }\tq^\tbeta;
\end{align*}
while the other half of the contribution $I_-$ is the following
\begin{align*}
I_-&= \sum_{\tbeta\in \tbL_\eff(L),w<0}
(-1)^{ (f+n) w} \frac{\prod_{m=1}^{-w-1}((f+n)(-w) + m)}{ (-w) \cdot(-w)!} \cdot\frac{1}{\prod_{i\in I\backslash (I_0\cup \{i_4\})} \langle D_i^*,\beta \rangle !} \cdot \frac{\prod_{m=-\infty}^{-w} m} {\prod_{m=-\infty}^{\langle D_{i_4}^*,\beta \rangle} m} 
\\ &\cdot  \frac{\prod_{m=-\infty}^{-w}(w+m)} {\prod_{m=-\infty}^{\langle D_{i_1}^*,\beta \rangle }(w+ m)}    \cdot \frac{\prod_{m=-\infty}^{0}((f+n+1)w+m)} {\prod_{m=-\infty}^{\langle D_{i_2}^*,\beta\rangle } (fw + m)} 
 \cdot \frac{\prod_{m=-\infty}^{0} (-(f+n)w + m)} {\prod_{m=-\infty}^{\langle D_{i_3}^*, \beta \rangle} ((-f-1)w + m)} q_0^w\prod_{a=1}^k q^\beta \\
&=\sum_{\tbeta\in \tbL_\eff(L), w<0}  \frac{(-1)^{\langle D_{i_3}^*+D_{r+3}^*,\tbeta\rangle}}{\prod_{i\in \tI\backslash \{i_2,i_3\}} \langle D_i^*,\tbeta \rangle!} \frac{\prod_{m=-\infty}^{ -\langle D_{i_3}^*,\tbeta \rangle -1} m }{\prod_{m=-\infty}^{ \langle D_{i_2}^*, \tbeta \rangle} m }\tq^\tbeta
\end{align*}
Therefore,
$$
F(Q;f) =I_++I_-=\sum_{\tbeta\in \tbL_\eff(L)}  \frac{(-1)^{\langle D_{i_3}^*+D_{r+3}^*,\tbeta\rangle}}{\prod_{i\in \tI\backslash \{i_2,i_3\}} \langle D_i^*,\tbeta \rangle!} \frac{\prod_{m=-\infty}^{ -\langle D_{i_3}^*,\tbeta \rangle -1} m }{\prod_{m=-\infty}^{ \langle D_{i_2}^*, \tbeta \rangle} m }\tq^\tbeta.
$$
It agrees with Equation \eqref{eqn:superpotential}.  This completes the proof of Conjecture  \ref{AVconj} when $L$ is a framed inner brane.

\section{Explicit Mirror Formulae}
\label{sec:formula}

In this section we list explicit mirror formulae for several examples, whose charge vectors are listed in Table \ref{table:charge}. The open and closed mirror maps are
given by
$$
Q_a= q_a\exp( \sum_{i\in I} l^{(a)}_i A_i (q)),\quad a=0,\ldots, k.
$$
The formulae for $$A_i(q)=\sum_{d_1,\ldots,d_k\in \bZ} a^i_{d_1,\ldots,d_k} \prod_{a=1}^k q_a^{d_a}.$$ are listed in Table \ref{table:A} at the end of this section. In particular, there is no mirror correction when $X$ is the resolved conifold. Notice $A_i(q)=0$ for $i\neq 1$ for $X$ other than 
the toric crepant resolution $Y_m$ of $(\cO_{\bP^1}(-1)\oplus \cO_{\bP^1}(-1))/\bZ_m$.

\begin{table}[h]
\begin{tabular}{|c|c|} \hline
$S$ &  $l^{(a)}$,  $a=1,\ldots,k$ \\ \hline \hline
$\cO_{\bP^1}(-1)$  & $l^{(1)}=(-1,-1,1,1)$\\ \hline
$\bP^2$ & $l^{(1)}= (-3,1,1,1)$ \\ \hline
$\bP^1\times \bP^1$ & $\begin{array}{l} l^{(1)}=(-2,1,1,0,0) \\ l^{(2)}=(-2,0,0,1,1)\end{array}$ \\ \hline
$dP_1$  & $\begin{array}{l} l^{(1)}=(-2,1,1,0,0) \\  l^{(2)}=(-1,0,-1,1,1)\end{array}$ \\ \hline
$dP_2$ &   $\begin{array}{l}    l^{(1)}=(-2,1,1,0,0,0)\\ l^{(2)}=(-2,0,0,1,1,0)\\ l^{(3)}=(-3,1,0,1,0,1) \end{array}$ \\ \hline
$dP_3$ & $\begin{array}{l}    l^{(1)}=(-2,1,1,0,0,0,0)\\ l^{(2)}=(-2,0,0,1,1,0,0)\\
l^{(3)}=(-3,1,0,1,0,1,0) \\ l^{(4)}=(-3,0,1,0,1,0,1) \end{array}$ \\ \hline
$Y_m$ & $\begin{array}{l} l^{(1)}=(1,1,0,-2,0,0,\dots,0)\\l^{(2)}=(0,0,1,-2,1,0,\dots, 0)\\l^{(3)}=(0,0,0,1,-2,1,\dots,0)\\ \dots \\l^{(m)}=(0,0,0,\dots,0,1,-2,1)\end{array}$ \\ \hline
\end{tabular}
\caption{charge vectors}
\label{table:charge}
\end{table}

\subsection{$\cO_{\bP^1}(-1)\oplus \cO_{\bP^1}(-1)$} \label{sec:conifold}
\subsubsection*{$B$-model:}
$\bT^\vee_q$ is
$$
x_1^{-1} x_2^{-1} x_3 x_4 = q=e^t.
$$
The mirror curves are listed in the following table.
\begin{center}
\begin{tabular}{|c|c|} \hline
Phase & Mirror curve \\\hline \hline
I & $\wy +1-q \wy -\wx \wy^{-f}=0$\\ \hline
II & $-\wx \wy^{-f}+ \wy+1 - q \wx \wy^{-f+1}=0$ \\ \hline
\end{tabular}
\end{center}

\subsubsection*{$A$-model:}
\begin{figure}[h]
\psfrag{X1=0}{\small $X_3=0$}
\psfrag{X2=0}{\small $X_4=0$}
\psfrag{X3=0}{\small $X_1=0$}
\psfrag{X4=0}{\small $X_2=0$}
\psfrag{c1}{\small $c_1$}
\psfrag{c2}{\small $c_2$}
\psfrag{r-c1}{\small $r-c_1$}
\includegraphics[scale=0.7]{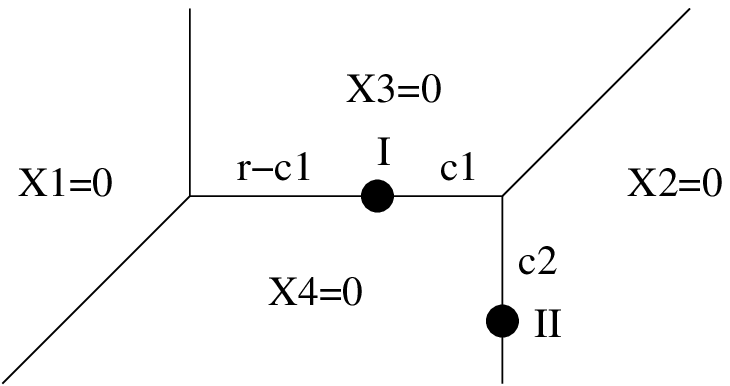}
\caption{Phase I: $r> c_1>0, c_2=0$. Phase II: $c_1=0, c_2>0$.}
\label{fig:conifold}
\end{figure}

$$
-|X_1|^2 - |X_2|^2 +|X_3|^2 +|X_4|^2 =r,
$$
where $r>0$.
$$
|X_4|^2 - |X_2|^2 = c_1,\quad
|X_1|^2 - |X_2|^2 = c_2. 
$$
Therefore
$$
(|X_1|^2, |X_2|^2, |X_3|^2, |X_4|^2)= (c_2, 0, r-c_1+c_2, c_1)+ s(1,1,1,1),\quad s\geq 0.
$$
See Figure \ref{fig:conifold}.

\subsubsection*{Mirror Formula:} 
$F(Q;f)=W(q;f)=\displaystyle{ \sum_{w\neq 0}\sum_d n_{d,w}(f) q_0^w q_1^d }$.
\begin{center}{\small
\begin{tabular}{|c|c|c|c|}\hline
Phase & $(i_1,i_2,i_3)$&   $n_{d,w}(f)$  & $\tbL_\eff(L)$  \\ \hline \hline
I & $(4,1,2)$ &  $\displaystyle{ (-1)^{fw+d}\frac{ \prod_{m=-d+1}^{d+w-1} (fw+m)}{w\cdot (w+d)!d!} }$ & $d, w+d\geq 0$ \\ \hline
II & $(1,2,4)$ & $\displaystyle{(-1)^{fw+d} \frac{\prod_{m=-d+1}^{-d+w-1}(fw+m)}{w\cdot(w-d)!d!} }$  & $w\geq d\geq 0$  \\ \hline
\end{tabular} }
\end{center}

\subsection{$K_{\bP^2}$}
\subsubsection*{$B$-model:}
$\bT_q^\vee$ is
$$
x_1^{-3} x_2 x_3 x_4= q=e^t.
$$
The mirror curves are listed in the following table.
\begin{center}
\begin{tabular}{|c|c|} \hline
Phase & Mirror curve \\\hline \hline
I & $1-\wx \wy^{-f}+\wy -q \wx^{-1} \wy^{f-1}=0$\\ \hline
II & $\wy+1-\wx \wy^{-f}-q\wx^{-1}\wy^{f+3}=0$ \\ \hline
III & $-\wx\wy^{-f}+\wy+1-q \wx^3\wy^{-3f-1}=0$ \\ \hline
\end{tabular}
\end{center}

\subsubsection*{$A$-model:}
\begin{figure}[h]
\psfrag{X0=0}{\small $X_1=0$}
\psfrag{X1=0}{\small $X_2=0$}
\psfrag{X2=0}{\small $X_3=0$}
\psfrag{X3=0}{\small $X_4=0$}
\psfrag{c1}{\small $c_1$}
\psfrag{c2}{\small $c_2$}
\psfrag{r-c1}{\small $r-c_1$}
\psfrag{r-c2}{\small $r-c_2$}
\psfrag{-c1}{\small $-c_1$}
\includegraphics[scale=0.7]{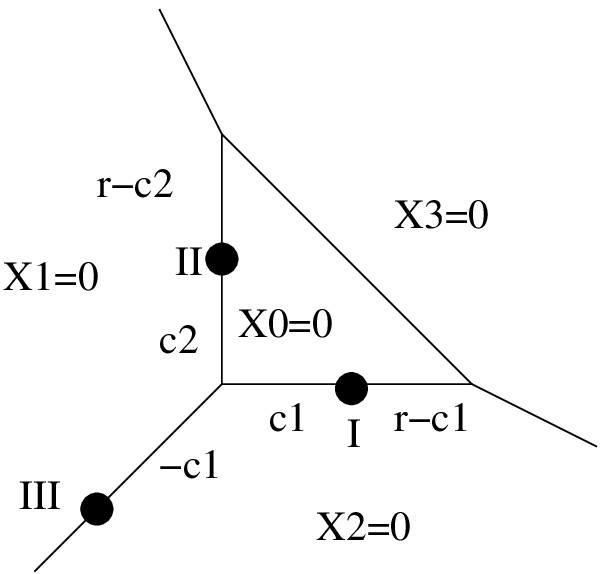}
\caption{Phase I: $r> c_1>0, c_2=0$. Phase II: $c_1=0, r> c_2>0$. Phase III: $c_1= c_2<0$.}
\label{fig:P2}
\end{figure}

$$
|X_2|^2 + |X_3|^2 +|X_4|^2 -3|X_1|^2 =r,
$$
where $r>0$.
$$
|X_2|^2 - |X_1|^2 = c_1,\quad
|X_3|^2 - |X_1|^2 = c_2. 
$$
Therefore
$$
(|X_1|^2, |X_2|^2, |X_3|^2, |X_4|^2)= (0, c_1, c_2, r-(c_1+c_2))+ s(1,1,1,1),\quad s\geq 0.
$$
See Figure \ref{fig:P2}.

\subsubsection*{Mirror Formula:} 
$F(Q;f)=W(q;f)=\displaystyle{ \sum_{w\neq 0}\sum_d n_{d,w}(f) q_0^w q_1^d }$.

\begin{center}{\small
\begin{tabular}{|c|c|c|c|}\hline
Phase & $(i_1,i_2,i_3)$ &  $n_{d,w}(f)$  & $\tbL_\eff(L)$ \\ \hline \hline
I   & $(2,3,1)$ & $\displaystyle{ (-1)^{fw+d}\frac{ \prod_{m=d+1}^{3d+w-1} (fw+m)}{w\cdot (w+d)!d!} }$ & $d, w+d\geq 0$ \\ \hline
II  & $(3,1,2)$ & $\displaystyle{(-1)^{fw+d} \frac{\prod_{m=-3d+1}^{-d+w-1}(fw+m)}{w\cdot(w+d)!d!} }$  & $d, w+d \geq 0$ \\ \hline
III & $(1,2,3)$ & $\displaystyle{(-1)^{fw+d} \frac{\prod_{m=d+1}^{-d+w-1}(fw+m)}{w\cdot(w-3d)!d!} }$  & $d \geq 0$, $w\geq 3d$ \\ \hline
\end{tabular} }
\end{center}

\bigskip

In Phase I, when $f=0$ we have
$$
\sum_{w\neq 0, d\geq 0, w+d\geq 0}N_{0,d\alpha +wb,w}(f) Q_0^w Q_1^d
= \sum_{w\neq 0, d\geq 0, w+d \geq 0} (-1)^{d} \frac{(w+3d-1)!}{w\cdot (w+d)! (d!)^2}  q_0^w q_1^d ,
$$
which agrees with the formula in \cite[Section 4.3]{GZ}.

In Phase III, when $f=0$ we have
$$
\sum_{w>0, d\geq 0} N_{0,d\alpha + wb,w}(f) Q_0^w Q_1^d
= \sum_{w> 0,d\geq 0, w\geq 3d}  (-1)^d \frac{(w-d-1)!}{w \cdot (d!)^2 (w-3d)!} q_0^w q_1^d,
$$
which agrees with \cite[Equation A.3]{LM}.

\subsection{$K_{\bP^1\times  \bP^1}=K_{\bF_0}$}
\subsubsection*{$B$-model:}
The torus $\bT_q^\vee$ is given by the following equations
\begin{align*}
x_1^{-2} x_2 x_3 = q_1=e^{t_1},\\
x_1^{-2} x_4 x_5 = q_2=e^{t_2}. 
\end{align*}
\begin{center}
The mirror curves are listed in the following table.
\begin{tabular}{|c|c|} \hline
Phase & Mirror curve \\\hline \hline
I & $\wy-q_1 \wx^{-1} \wy^{f+2}-\wx \wy^{-f}+1+q_2\wy^2=0$\\ \hline
II & $1+q_1 \wy^{-1} +\wy -\wx \wy^{-f}-q_2 \wx^{-1}\wy^f=0$ \\ \hline
III & $-\wx \wy^{-f}+q_1 \wx^2\wy^{-2f}+1+\wy+q_2 \wx^2 \wy^{-2f-1}=0$ \\ \hline
\end{tabular}
\end{center}

\subsubsection*{$A$-model:}
\begin{figure}[h]
\psfrag{X0=0}{\small $X_1=0$}
\psfrag{X1=0}{\small $X_2=0$}
\psfrag{X2=0}{\small $X_3=0$}
\psfrag{X3=0}{\small $X_4=0$}
\psfrag{X4=0}{\small $X_5=0$}
\psfrag{c1}{\small $c_1$}
\psfrag{c2}{\small $c_2$}
\psfrag{r1-c1}{\small $r_1-c_1$}
\psfrag{r2-c2}{\small $r_2-c_2$}
\psfrag{-c1}{\small $-c_1$}
\includegraphics[scale=0.7]{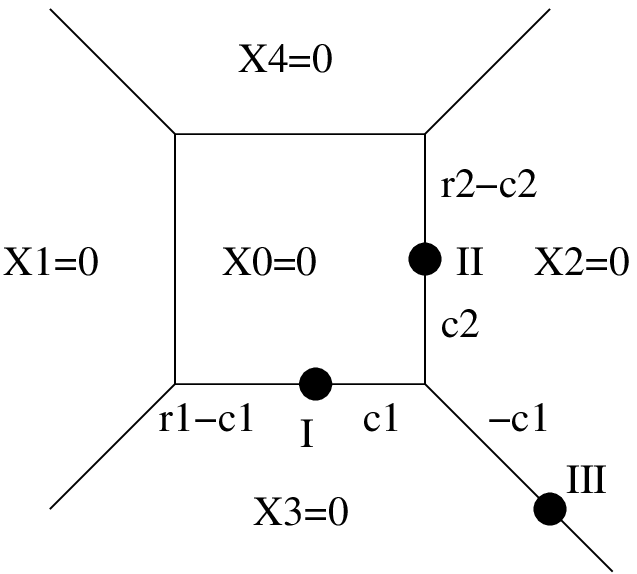}
\caption{Phase I: $r_1> c_1>0, c_2=0$. Phase II: $c_1=0, r_2> c_2>0$. Phase III: $c_1= c_2<0$}
\label{fig:F0}
\end{figure}

$$
|X_2|^2 + |X_3|^2 -2|X_1|^2 =r_1,\quad |X_4|^2 + |X_5|^2-2|X_1|^2 = r_2
$$
where $r_1, r_2 >0$.
$$
|X_3|^2 - |X_1|^2 = c_1,\quad
|X_4|^2 - |X_1|^2 = c_2. 
$$
Therefore
$$
(|X_1|^2, |X_2|^2,\ldots, |X_5|^2)= (0,r_1-c_1 , c_1, c_2,r_2-c_2)+ s(1,1,1,1,1),\quad s\geq 0.
$$
See Figure \ref{fig:F0}.

\subsubsection*{Mirror Formula:} 
$F(Q;f)=W(q;f)=\displaystyle{ \sum_{w\neq 0}\sum_{d_1,d_2} n_{d_1,d_2,w}(f) q_0^w q_1^{d_1} q_2^{d_2} }$.
\begin{center} {\small
\begin{tabular}{|c|c|c|c|}\hline
Phase & $(i_1,i_2,i_3)$ &  $n_{d_1,d_2,w}(f)$  & $\tbL_\eff(L)$ \\ \hline \hline
I & $(3,1,4)$ & $\displaystyle{ (-1)^{fw+d_2}\frac{ \prod_{m=-2d_1-2d_2+1}^{-d_2+w-1} (fw+m)}{w\cdot (w+d_1)!d_1!d_2! } }$ & $d_1,d_2, w+d_1\geq 0$ \\ \hline
II & $(4,3,1)$& $\displaystyle{(-1)^{fw} \frac{\prod_{m=d_1+1}^{2d_1+2d_2+w-1}(fw+m)}{w\cdot(w+d_2)!d_1!d_2!} }$  & $d_1,d_2,w+d_2 \geq 0$ \\ \hline
III & $(1,4,3)$ & $\displaystyle{(-1)^{fw+d_1} \frac{\prod_{m=d_2+1}^{-d_1+w-1}(fw+m)}{w\cdot(w-2d_1-2d_2)!d_1!d_2!} }$  
& $d_1,d_2 \geq 0$, $w\geq 2d_1+2d_2$ \\ \hline
\end{tabular} }
\end{center}

\subsection{$K_{dP_1}= K_{\bF_1}$}
\subsubsection*{$B$-model:}
The torus $\bT_q^\vee$ is given by the following equations
\begin{eqnarray*}
&x_1^{-2} x_2 x_3 = q_1=e^{t_1},\\
&x_1^{-1} x_3^{-1} x_4 x_5 = q_2=e^{t_2}.
\end{eqnarray*}
The mirror curves are listed in the following table.
\begin{center}
\begin{tabular}{|c|c|} \hline
Phase & Mirror curve \\\hline \hline
I & $\wy+1+q_1 \wy^2-\wx \wy^{-f}-q_1 q_2 \wx^{-1}\wy^{f+3}=0$\\ \hline
II & $1-\wx \wy^{-f}-q_1 \wx^{-1} \wy^f+\wy-q_1 q_2 \wx^{-1}\wy^{f-1}=0$ \\ \hline
III & $1+q_1 \wy^{-1}+\wy-\wx \wy^{-f}- q_2 \wx^{-1}\wy^{f+1}=0$ \\ \hline
IV & $-\wx \wy^{-f}+\wy+q_1 \wx^2 \wy^{-2f-1}+1 - q_1 q_2 \wx^3 \wy^{-3f-1}=0$ \\ \hline
V & $-\wx \wy^{-f}+q_1 \wx^2 \wy^{-2f}+1+\wy- q_2 \wx \wy^{-f-1}=0$ \\ \hline
\end{tabular}
\end{center}

\subsubsection*{$A$-model:}
\begin{figure}[h]
\psfrag{X0=0}{\small $X_1=0$}
\psfrag{X1=0}{\small $X_4=0$}
\psfrag{X2=0}{\small $X_5=0$}
\psfrag{X3=0}{\small $X_3=0$}
\psfrag{X4=0}{\small $X_2=0$}
\psfrag{c1}{\small $c_1$}
\psfrag{c2}{\small $c_2$}
\psfrag{r1-c2}{\small $r_1-c_2$}
\psfrag{r2-c1}{\small $r_2-c_1$}
\psfrag{r1+r2-c1}{\small $r_1+r_2-c_1$}
\includegraphics[scale=0.7]{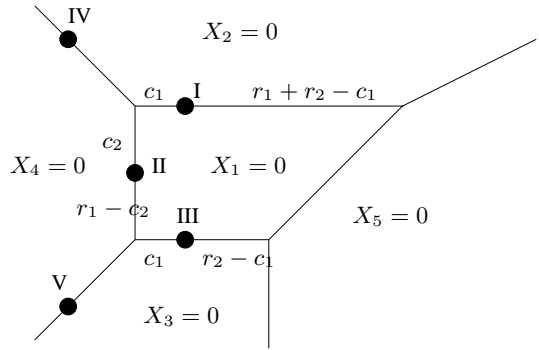}
\caption{Phase I: $r_1 + r_2 > c_1>0, c_2=0$. Phase II: $c_1=0, r_1 > c_2>0$. Phase III: $r_2 >c_1>0,  c_2=r_1$.
Phase IV: $c_1=c_2<0$. Phase V: $c_1=r_1-c_2<0$.}
\label{fig:F1}
\end{figure}

$$
|X_2|^2 + |X_3|^2 -2|X_1|^2 =r_1,\quad |X_4|^2 + |X_5|^2-|X_1|^2 -|X_3|^2 = r_2
$$
where $r_1, r_2>0$.
$$
|X_4|^2 - |X_1|^2 = c_1,\quad
|X_2|^2 - |X_1|^2 = c_2. 
$$
Therefore
$$
(|X_1|^2, |X_2|^2, \ldots, |X_5|^2)= (0, c_2, r_1-c_2,c_1  ,r_1+r_2-c_1-c_2)+ s(1, 1,1,1,1),\quad s\geq 0.
$$
See Figure \ref{fig:F1}.

\subsubsection*{Mirror Formula:} 
$F(Q;f)=W(q;f)=\displaystyle{ \sum_{w\neq 0}\sum_{d_1,d_2} n_{d_1,d_2,w}(f) q_0^w q_1^{d_1} q_2^{d_2} }$.

\begin{center}{\small
\begin{tabular}{|c|c|c|c|}\hline
Phase & $(i_1,i_2,i_3)$ &  $n_{d_1,d_2,w}(f)$  & $\tbL_\eff(L)$ \\ \hline \hline
I & $(4,1,2)$ & $\displaystyle{ (-1)^{fw+d_1}\frac{ \prod_{m=-2d_1-d_2+1}^{-d_1+w-1} (fw+m)}{w\cdot (w+d_2)!(d_1-d_2)!d_2! } }$ & $w+d_2\geq 0$, $d_1\geq d_2\geq 0$ \\ \hline
II & $(2,4,1)$& $\displaystyle{(-1)^{fw+d_2} \frac{\prod_{m=d_2+1}^{2d_1+d_2+w-1}(fw+m)}{w\cdot(w+d_1)!(d_1-d_2)!d_2!} }$  & $w+d_1\geq 0$, $d_1\geq d_2\geq 0$ \\ \hline
III & $(4,3,1)$ & $\displaystyle{(-1)^{fw+d_2} \frac{\prod_{m=d_1-d_2+1}^{2d_1+d_2+w-1}(fw+m)}{w\cdot(w+d_2)!d_1!d_2!} }$  
& $d_1,d_2 \geq 0$, $w+d_2\geq 0$ \\ \hline
IV & $(1,2,4)$& $\displaystyle{(-1)^{fw+d_2} \frac{\prod_{m=d_1+1}^{-d_2+w-1}(fw+m)}{w\cdot(w-2d_1-d_2)!(d_1-d_2)!d_2!} }$  & $d_1\geq d_2\geq 0$, $w\geq 2d_2+d_2$ \\ \hline
V & $(1,4,3)$ & $\displaystyle{(-1)^{fw+d_1+d_2} \frac{\prod_{m=d_2+1}^{-d_1+d_2+w-1}(fw+m)}{w\cdot(w-2d_1-d_2)!d_1!d_2!} }$  
& $d_1,d_2 \geq 0$, $w\geq 2d_1+d_2$ \\ \hline
\end{tabular} }
\end{center}

\subsection{$K_{dP_2}$}
\subsubsection*{$B$-model:}
The torus $\bT_q^\vee$ is given by
\begin{eqnarray*}
&x_1^{-2} x_2 x_3 = q_1=e^{t_1},\\ &x_1^{-2}x_4x_5=q_2=e^{t_2},\\ &x_1^{-3}x_2x_4x_6=q_3=e^{t_3}.
\end{eqnarray*}
The mirror curves are listed in the following table.
\begin{center}
\begin{tabular}{|c|c|} \hline
Phase & Mirror curve \\\hline \hline
I & $\wy+1+q_1 \wy^{2}- q_2 \wx^{-1} \wy^{f+2}-\wx \wy^{-f}-(q_2)^{-1}q_3 \wx \wy^{-f+1}=0$\\ \hline
II & $1-\wx \wy^{-f}+q_1 \wx^2 \wy^{-2f}+ q_2 \wy^{-1}+\wy-(q_2)^{-1}q_3 \wx^{-1} \wy^{f+1}=0$ \\ \hline
III & $1-(q_2)^{-1} q_3 \wx \wy^{-f-1}- q_1 q_2 (q_3)^{-1} \wy^{f+1}-q_2 \wx^{-1} \wy^{f}-\wx \wy^{-f}+\wy=0$ \\ \hline
IV & $-\wx \wy^{-f}+1+q_1 \wx^2 \wy^{-2f}+\wy+ q_2 \wx^2 \wy^{-2f-1}-q_3 \wx^3 \wy^{-3f-1}=0$ \\ \hline
V & $-\wx \wy^{-f}+\wy+q_1 \wx^2 \wy^{-2f-1}+q_2 \wx^2 \wy^{-2f}+1-(q_2)^{-1} q_3 \wx \wy^{-f-1}=0$ \\ \hline
VI & $-xy^{-f}-q_2^{-1}q_3 x y^{-f+1} - q_1 q_2 q_3^{-1} x y^{-f-1} + q_2 x^2 y^{-2f-1} + y +1 =0 $ \\ \hline
\end{tabular}
\end{center}

\subsubsection*{$A$-model:}
\begin{eqnarray*}
|X_2|^2+|X_3|^2-2|X_1|^2 &=& r_1,\\
|X_4|^2+|X_5|^2-2|X_1|^2 &=& r_2,\\
|X_2|^2+|X_4|^2+|X_6|^2-3|X_1|^2 &=& r_3,
\end{eqnarray*}
where $r_1,r_2,r_3>0$, and $\max(r_1,r_2)<r_3<r_1+r_2$.
$$|X_2|^2-|X_1|^2=c_1,\quad |X_5|^2-|X_1|^2=c_2.$$
Therefore
$$
(|X_1|^2,|X_2|^2,\ldots,|X_6|^2)=(0,c_1,r_1-c_1 , r_2-c_2,c_2, r_3-r_2-c_1+c_2)+ t(1,1,1,1,1,1),\quad t\geq 0.
$$

\begin{figure}[h]
\psfrag{X0=0}{\small $X_1=0$}
\psfrag{X1=0}{\small $X_2=0$}
\psfrag{X2=0}{\small $X_3=0$}
\psfrag{X3=0}{\small $X_4=0$}
\psfrag{X4=0}{\small $X_5=0$}
\psfrag{X5=0}{\small $X_6=0$}
\psfrag{c1}{\small $c_1$}
\psfrag{c2}{\small $c_2$}
\psfrag{c2=c1-r3+r2}{\small $c_2$}
\psfrag{r2-c2=c1<0}{\small $-c_1$}
\psfrag{c1=c2<0}{\small $\quad \quad\quad-c_1$}
\psfrag{c1=r3-r2, c2<0}{\small $-c_2$}
\includegraphics[scale=0.7]{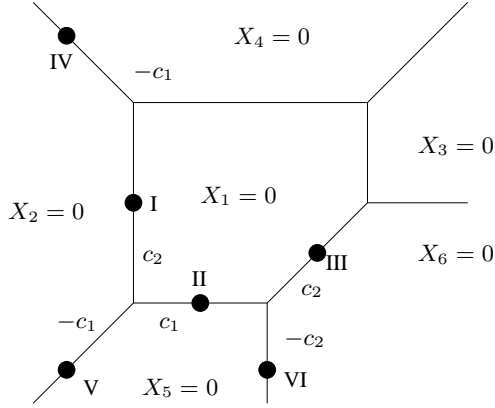}
\caption{Phase I: $r_2>c_2>0, c_1=0$. Phase II: $r_3-r_2>c_1>0, c_2=0$. Phase III: $r_3-r_2<c_1<r_1, c_1-c_2=r_3-r_2$. Phase IV: $r_2-c_2=c_1<0$ Phase V: $c_1=c_2<0$. Phase VI: $c_1=r_3-r_2, c_2<0$.}
\label{fig:dP2}
\end{figure}

See Figure \ref{fig:dP2}.

\subsubsection*{Mirror Formula:} 
$F(Q;f)=W(q;f)=\displaystyle{ \sum_{w\neq 0}\sum_{d_1,d_2, d_3} n_{d_1,d_2,d_3,w}(f) q_0^w \prod_{a=1}^3 q_a^{d_a} }$.

\begin{center} {\small
\begin{tabular}{|c|c|c|c|} \hline
Phase & $(i_1,i_2,i_3)$ & $n_{d_1,d_2,d_3,w}(f)$ & $\tbL_\eff(L)$ \\ \hline \hline
I & $(5,1,2)$ &
$\displaystyle{ \frac{(-1)^{fw+d_1+d_3} \prod_{m=-2d_1-2d_2-3d_3+1}^{-d_1- d_3 + w-1}(fw+m)}{w(w+d_2)!d_1!d_3!(d_2+d_3)!} }$ & 
$\begin{array}{c} d_1, d_3, d_2+d_3\geq 0 \\ w+d_2\geq 0 \end{array}$  \\ \hline
II & $(2,5,1)$ &
$\displaystyle{ \frac{(-1)^{fw+d_3}\prod_{m=d_2+1}^{2d_1+2d_2+3d_3 +w -1}(fw+m)}{w(w+d_1+d_3)!d_1!d_3!(d_2+d_3)!} }$ & 
$\begin{array}{c} d_1, d_3, d_2+d_3\geq 0\\ w+d_1+d_3\geq 0\end{array}$ \\ \hline
III & $(5,6,1)$ &
$\displaystyle{ \frac{(-1)^{fw+d_3}\prod_{m=d_3+1}^{2d_1+2d_2+3d_3 +w -1}(fw+m)}{w(w+d_2)!d_1!(d_1+d_3)!(d_2+d_3)!} }$ &
$\begin{array}{c} d_1, d_1+d_3, d_2+d_3\geq 0 \\ w+d_2\geq 0\end{array}$ \\ \hline
IV & $(1,4,2)$ &
$\displaystyle{ \frac{(-1)^{fw+d_1+d_3}\prod_{m=d_2+ d_3+1}^{-d_1-d_3 +w -1}(fw+m)}{w(w-2d_1-2d_2-3d_3)!d_1!d_2!d_3!} }$ & 
$\begin{array}{c} d_1,d_2,d_3\geq 0 \\ w\geq 2d_1+2d_2+3d_3\end{array}$  \\ \hline
V  & $(1,2,5)$ &
$\displaystyle{  \frac{(-1)^{fw+d_2} \prod_{m=d_1+ d_3+1}^{-d_2 +w -1}(fw+m)}{w(w-2d_1-2d_2-3d_3)!d_1!(d_2+d_3)!d_3!} }$ &
$\begin{array}{c} d_1, d_3, d_2+d_3\geq 0\\ w\geq 2d_1+2d_2+3d_3\end{array}$ \\ \hline
VI & $(1,5,6)$ &
$\displaystyle{ \frac{(-1)^{fw+d_3}\prod_{m=d_1+ d_3+1}^{-d_3 +w -1}(fw+m)}{w(w-2d_1-2d_2-3d_3)!d_1!(d_1+d_3)!(d_2+d_3)!} }$  &
$\begin{array}{c} d_1, d_1+d_3, d_2+d_3\geq 0 \\  w\geq 2d_1+2d_2+3d_3 \end{array}$ \\ \hline
\end{tabular} }
\end{center}

\subsection{$K_{dP_3}$} \label{sec:dP3}
\subsubsection*{$B$-model:}
The torus $\bT_q^\vee$ is given by
$$
x_1^{-2} x_2 x_3=q_1=e^{t_1},\quad x_1^{-2}x_4 x_5=q_2=e^{t_2},\quad x_1^{-3}x_2 x_4 x_6=q_3=e^{t_3}, \quad x_1^{-3} x_3 x_5 x_7=q_4=e^{t_4}.
$$
The mirror curves are listed in the following table.
\begin{center}
\begin{tabular}{|c|c|} \hline
Phase & Mirror curve \\\hline \hline
I & $1+\wy+q_1 \wy^{-1}-q_2 \wx^{-1} \wy^f-\wx \wy^{-f}-(q_2)^{-1}q_3 \wx \wy^{-f-1}-(q_1)^{-1}q_4 \wx^{-1} \wy^{f+1}=0$\\ \hline
II & $-\wx \wy^{-f}+1+q_1 \wx^2 \wy^{-2f}+q_2 \wx^2 \wy^{-2f-1}+\wy-(q_2)^{-1}q_3 \wx \wy^{-f+1}-(q_1)^{-1} q_4 \wx \wy^{-f-1}=0$ \\ \hline
\end{tabular}
\end{center}

\subsubsection*{$A$-model:}
\begin{eqnarray*}
|X_2|^2+|X_3|^2-2|X_1|^2&=&r_1,\\
|X_4|^2+|X_5|^2-2|X_1|^2&=&r_2,\\
|X_2|^2+|X_4|^2+|X_6|^2-3|X_1|^2&=&r_3,\\
|X_3|^2+|X_5|^2+|X_7|^2-3|X_1|^2&=&r_4,
\end{eqnarray*}
where $r_1,r_2>0$, $\max(r_1,r_2)<r_3<r_1+r_2$, $\max(r_1,r_2)<r_4<r_1+r_2$.
$$|X_2|^2-|X_1|^2=c_1,\quad |X_5|^2-|X_1|^2=c_2.$$ 
Therefore
\begin{eqnarray*}
&(|X_1|^2,|X_2|^2,|X_3|^2,|X_4|^2,|X_5|^2,|X_6|^2, |X_7|^2)\\=&(0,c_1,r_1-c_1, r_2-c_2,c_2, r_3-r_2-c_1+c_2,r_4-r_1+c_1-c_2)+t(1,1,1,1,1,1),\quad t\geq 0.
\end{eqnarray*}

\begin{figure}[h]
\psfrag{X0=0}{\small $X_1=0$}
\psfrag{X1=0}{\small $X_2=0$}
\psfrag{X2=0}{\small $X_3=0$}
\psfrag{X3=0}{\small $X_4=0$}
\psfrag{X4=0}{\small $X_5=0$}
\psfrag{X5=0}{\small $X_6=0$}
\psfrag{X6=0}{\small $X_7=0$}
\psfrag{c2}{\small $c_1$}
\psfrag{c1=c2<0}{\small $\quad\quad \quad -c_1$}
\includegraphics[scale=0.7]{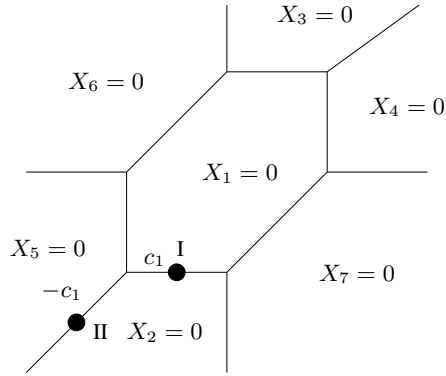}
\caption{Phase I: $0<c_1<r_4-r_1$. Phase II: $c_1= c_2>0$.}
\label{fig:dP3}
\end{figure}

See Figure \ref{fig:dP3}.

\subsubsection*{Mirror Formula:} 
$F(Q;f)=W(q;f)=\displaystyle{ \sum_{w\neq 0}\sum_{d_1,\ldots, d_4} n_{d_1,d_2,d_3,d_4,w}(f) q_0^w \prod_{a=1}^4 q_a^{d_a} }$.

\begin{center} {\small
\begin{tabular}{|c|c|c|c|} \hline
Phase & $(i_1,i_2,i_3)$ & $n_{d_1,d_2,d_3,d_4,w}(f)$ & $\tbL_\eff(L)$  \\ \hline \hline
I & $(5,2,1)$ & {\footnotesize $\displaystyle{ \frac{(-1)^{fw+d_3+d_4}\prod_{m=d_1+d_3+1}^{2 d_1+2d_2+3d_3+3d_4 +w -1}(fw+m)}{w(w+d_2+d_4)!d_3!d_4!(d_2+d_3)! (d_1+d_4)!} } $} 
& $\begin{array}{c} d_3,d_4, d_1+d_4, d_2+d_3\geq 0\\ w+d_2+d_4\geq 0\end{array}$ \\  \hline
II & $(1,5,2)$ &  {\footnotesize $\displaystyle{\frac{(-1)^{fw+d_2+d_4}\prod_{m=d_2+d_4+1}^{-d_1-d_3+w -1}(fw+m)}
{w(w-2d_1-2d_2-3 d_3 -3 d_4) )!d_3!d_4!(d_2+d_3)! (d_1+d_4)!} }$ }  
& $\begin{array}{c} d_3,d_4, d_1+d_4, d_2+d_3\geq 0\\ w\geq 2d_1+2d_2+3d_3+3d_4 \end{array}$ \\ \hline
\end{tabular} }
\end{center}

\subsection{The toric crepant resolution $Y_m$ of $(\cO_{\bP^1}(-1)\oplus \cO_{\bP^1}(-1))/\bZ_m$} \label{sec:Xm}

\subsubsection*{$B$-model:}
The torus $\bT_q^\vee$ is given by
\begin{align*}
x_1 x_2 x_4^{-2}&=q_1,\\
x_3 x_5 x_4^{-2}&=q_2,\\
x_4 x_6 x_5^{-2}&=q_3,\\
&\dots \\
x_{m+1} x_{m+3} x_{m+2}^{-2}&=q_m.
\end{align*}
One can write down the equations of the framed mirror curves \eqref{eqn:mirror-curve} prescribed in Section \ref{sec:hori-vafa-mirror}. We do not list them here because they are too long.

\subsubsection*{$A$-model:}
\begin{eqnarray*}
|X_1|^2+|X_2|^2-2|X_4|^2=r_1,\\
|X_3|^2+|X_5|^2-2|X_4|^2=r_2,\\
|X_4|^2+|X_6|^2-2|X_5|^2=r_3,\\
|X_5|^2+|X_7|^2-2|X_6|^2=r_4,\\
\dots,\\
|X_{m+1}|^2+|X_{m+3}|^2-2|X_{m+2}|^2=r_m,
\end{eqnarray*}
where $r_1,\dots, r_m >0$.

The Aganagic-Vafa A-branes in each phase are given as follows:
\begin{itemize}
\item $\rI_0$:  $$|X_1|^2-|X_4|^2=|X_3|^2- |X_4|^2=-c,\quad c>0$$
\item $\rI_b$, $1\le b \le m-1$: $$|X_1|^2 - |X_{b+3}|^2=0,\ |X_{b+2}|^2-|X_{b+3}|^2=c, \quad 0<c<r_{b+1}.$$
\item $\rI_m$: $$|X_1|^2 - |X_{m+3}|^2=0,\ |X_{m+2}|^2-|X_{m+3}|^2=c,\quad c>0.$$
\item $\rII_b$, $b=1,2$: $$|X_1|^2- |X_{b+3}|^2=c, \ |X_{b+2}|^2-|X_{b+3}|^2=0,\quad 0<c<r_1.$$
\item $\rII_b$, $2\le b \le m$: \begin{eqnarray*}&|X_1|^2- |X_{b+3}|^2=c, \ |X_{b+2}|^2-|X_{b+3}|^2=0,\\ & 0<c<r_1+2(r_3+2r_4+\dots + (b-2)r_b).\end{eqnarray*}
\end{itemize}

\begin{figure}[h]
\psfrag{X1}{\tiny $X_1=0$} 
\psfrag{X2}{\tiny $X_2=0$}
\psfrag{X3}{\tiny $X_3=0$} 
\psfrag{X4}{\tiny $X_4=0$}
\psfrag{X5}{\tiny $X_5=0$}
\psfrag{X6}{\tiny $X_6=0$}
\psfrag{X7}{\tiny $X_7=0$}
\psfrag{X8}{\tiny $X_8=0$}
\psfrag{r1}{\tiny $r_1$}
\psfrag{r2}{\tiny $r_2$}
\psfrag{r3}{\tiny $r_3$}
\psfrag{r4}{\tiny $r_4$}
\psfrag{r5}{\tiny $r_5$}
\psfrag{I0}{\tiny $\rI_0$}b
\psfrag{I1}{\tiny $\rI_1$}
\psfrag{I2}{\tiny $\rI_2$}
\psfrag{I3}{\tiny $\rI_3$}
\psfrag{I4}{\tiny $\rI_4$}
\psfrag{I5}{\tiny $\rI_5$}
\psfrag{II1}{\tiny $\rII_1$}
\psfrag{II2}{\tiny $\rII_2$}
\psfrag{II3}{\tiny $\rII_3$}
\psfrag{II4}{\tiny $\rII_4$}
\psfrag{II5}{\tiny $\rII_5$}
\includegraphics[scale=0.35]{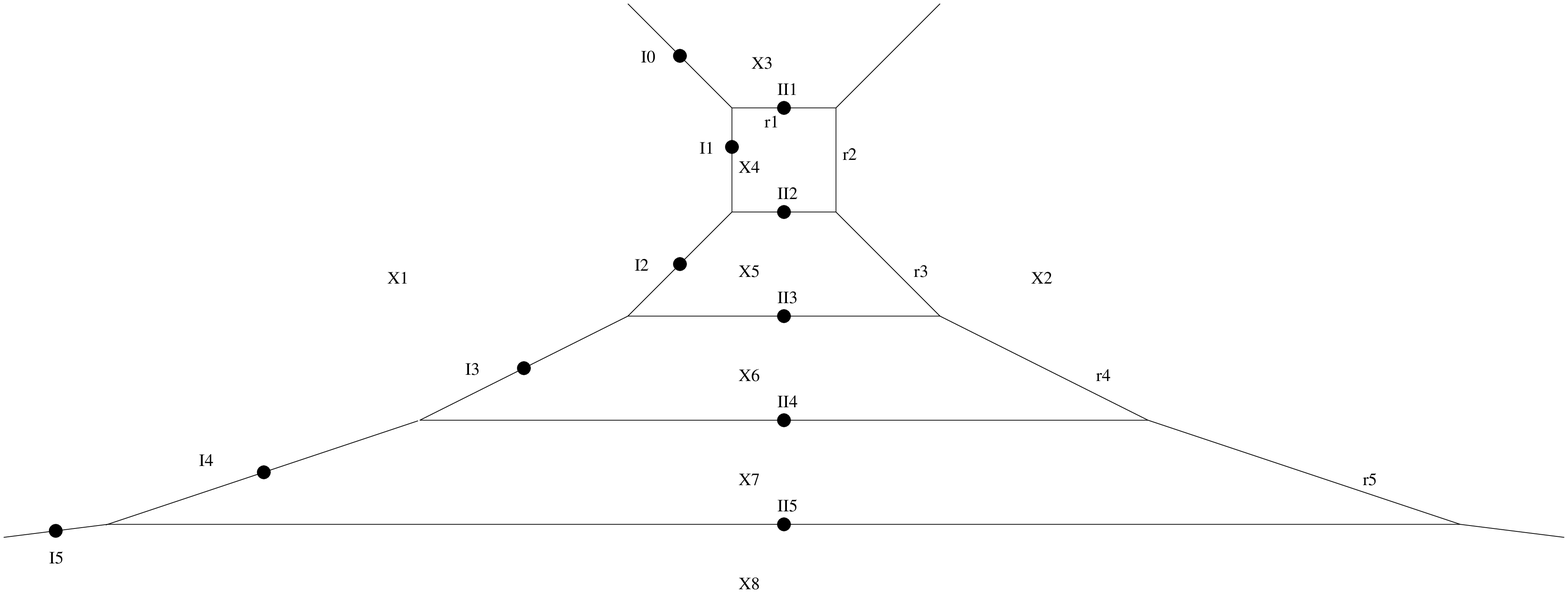}
\caption{$\mu'(X^1_5)=\mu'((\cO(-1)\oplus \cO(-1))/\bZ_5)^1)$.}
\label{fig:resolution}
\end{figure}

As an example, the image of the $1$-skeleton $X^1_5$ under the moment map $\mu': X\to (\ft'_\bR)^*$ is in Figure \ref{fig:resolution}.

\newpage

\subsubsection*{Mirror Formula:} 
$F(Q;f)=W(q;f)=\displaystyle{ \sum_{w\neq 0}\sum_{d_1,\ldots, d_m} n_{d_1,\dots,d_m,w}(f) q_0^w \prod_{a=1}^m q_a^{d_a} }$.

The mirror formulae for some phases are listed in the following table.
\begin{center} {\small
\begin{tabular}{|c|c|c|c|} \hline
Phase & $(i_1,i_2,i_3)$ & $n_{d_1,\dots, d_m,w}(f)$ & $\tbL_\eff(L)$  \\ \hline \hline
$\rI_0$ & $(4,3,1)$ & 
{\footnotesize $\begin{array}{l} (-1)^{fw+d_1}\prod_{m=d_2+1}^{-d_1+w-1} (fw+m) \\ 
\cdot \displaystyle{ \frac{1}{w(w-2d_1-2d_2+d_3)!d_1!} } \\ 
\cdot \displaystyle{ \frac{1}{(d_2-2d_3+d_4)!(d_3-2d_4+d_5)!} }\\  
\quad \dots \\ 
\cdot \displaystyle{ \frac{1}{(d_{m-2}-2d_{m-1}+d_m)!} \cdot \frac{1}{(d_{m-1}-2d_m)! d_m!} }  \end{array}$} 
& {\footnotesize $\begin{array}{l} w\geq 2d_1+2d_2-d_3,\\ d_1\geq 0,\\ d_2-2d_3+d_4\geq 0,\\ d_3-2d_4+d_5\geq 0,\\ \dots,\\ d_{m-2}-2d_{m-1}+d_m \geq 0,\\ d_{m-1}-2d_m\geq 0,\\ d_m\geq 0. \end{array}$}
\\  \hline
$\begin{array}{c}\rI_b\\ 3\leq b \leq m-2 \end{array}$ & $(b+2,1,b+3)$ &  {\footnotesize $\begin{array}{l} (-1)^{fw+d_b-2d_{b+1}+d_{b+2}}\\ 
\cdot \prod_{m=d_1+1}^{-d_b+2d_{b+1}-d_{b+2}+w-1} (fw+m) \\ 
\cdot \displaystyle{ \frac{1}{w(w+d_{b-1}-2d_b+d_{b+1})!d_1!d_2!} } \\ 
\cdot \displaystyle{ \frac{1}{(-2d_1-2d_2+d_3)!(d_3-2d_4+d_5)!} } \\
\quad \dots \\ 
\cdot \displaystyle{ \frac{1}{(d_{b-2}-2d_{b-1}+d_b)! (d_{b+1}-2d_{b+2}+d_{b+3})!} }\\ 
\quad \dots\\ 
\displaystyle{ \frac{1}{(d_{m-2}-2d_{m-1}+d_m)! (d_{m-1}-2d_m)! d_m!} }  \end{array} $  }  
& {\footnotesize $\begin{array}{l}\\ w+d_{b-1}-2d_b+d_{b+1}\geq 0,\\ d_1, d_2 \geq 0,\\-2d_1-2d_2+d_3\geq 0, \\ \dots \\ d_{b-2}-2d_{b-1}+d_b\geq 0,\\ d_{b+1}-2d_{b+2}+d_{b+3}\geq 0,\\ \dots,\\ d_{m-2}-2d_{m-1}+d_m\geq 0, \\ d_{m-1}-2d_m\geq 0,\\ d_m\geq 0 \\ { } \\ { } \end{array}$ }\\ \hline
$\begin{array}{c} \rII_b\\3\leq b \leq m-2 \end{array}$ & $(1,b+3,b+2)$ &  {\footnotesize $\begin{array}{l} (-1)^{fw+d_{b-1}-2d_{b}+d_{b+1}}\\ 
\cdot \prod_{m=d_{b+1}-2d_{b+2}+d_{b+3}+1}^{-d_{b-1}+2d_{b}-d_{b+1}+w-1} (fw+m)\\ 
\cdot \displaystyle{ \frac{1}{w(w+d_{1})!d_1!d_2!} } \\ 
\cdot \displaystyle{ \frac{1}{(-2d_1-2d_2+d_3)!(d_3-2d_4+d_5)!} }  \\ 
\quad \dots \\ 
\cdot \displaystyle{ \frac{1}{(d_{b-2}-2d_{b-1}+d_b)! (d_{b+1}-2d_{b+2}+d_{b+3})!}} \\ 
\quad \dots\\ 
\cdot \displaystyle{ \frac{1}{(d_{m-2}-2d_{m-1}+d_m)! (d_{m-1}-2d_m)! d_m!} } \end{array} $ }    &  
{\footnotesize $\begin{array}{l} \\ w+d_1 \geq 0,\\ d_1, d_2 \geq 0,\\-2d_1-2d_2+d_3\geq 0, \\ \dots \\ d_{b-2}-2d_{b-1}+d_b\geq 0,\\ d_{b+1}-2d_{b+2}+d_{b+3},\\ \dots,\\ d_{m-2}-2d_{m-1}+d_m\geq 0, \\ d_{m-1}-2d_m\geq 0,\\ d_m\geq 0 \\ { } \\ { } \end{array}$ }
\\ \hline
\end{tabular} }
\end{center}

\newpage

\begin{center}{\small 
\begin{table}[h]
\begin{tabular}{|c|c|c|c|}

\hline
$S$ & $i$ &  $a^i_{d_1,\ldots,d_k}$  &  constraints on $d_a$ in the sum  \\ \hline \hline

$\cO_{\bP^1}(-1)$ & $i=1,\dots, 4$ & $0$ &   \\ \hline

\multirow{2}{*}{$\bP^2$} & $1$ & $\displaystyle{\frac{(-1)^{d_1-1}(3d_1-1)!}{ (d_1!)^3}  }$ & $d_1>0$ \\ \cline{2-4}
& $2,3,4$ & 0 & \\ \hline

\multirow{2}{*}{$\bP^1\times \bP^1$} & 1 &$\displaystyle{ \frac{-(2d_1+2d_2-1)!}{(d_1!)^2(d_2!)^2}  }$ & 
$\begin{array}{c} d_1, d_2\geq 0 \\ (d_1,d_2)\neq (0,0)\end{array}$  \\ \cline{2-4}
& $2,\cdots, 5$ & 0 & \\ \hline

\multirow{2}{*}{$dP_1$} & $1$ & $\displaystyle{ \frac{(-1)^{d_2-1}(2d_1+d_2-1)!}{d_1!(d_1-d_2)!(d_2!)^2}  }$ &
$\begin{array}{c} d_1\geq d_2\geq 0 \\ (d_1,d_2)\neq (0,0)\end{array} $ \\ \cline{2-4}
& $2,\cdots, 5$ & 0 & \\ \hline

\multirow{2}{*}{$dP_2$} & $1$ &  $\displaystyle{ \frac{(-1)^{d_3-1}(2d_1+2d_2+3d_3-1)!}{d_1!d_2!d_3!(d_1+d_3)!(d_2+d_3)!} 
}$  & $\begin{array}{c} d_1,d_2,d_3\geq 0\\(d_1,d_2,d_3)\neq (0,0,0)\end{array}$    \\ \cline{2-4}
& $2,\cdots, 6$ & 0 & \\ \hline

\multirow{2}{*}{$dP_3$} & $1$ & $\displaystyle{  \frac{(-1)^{d_3+d_4 -1}(2d_1+2d_2+3d_3+3d_4 -1)!}
{(d_1+d_3)!(d_1+d_4)!(d_2+d_3)!(d_2+d_4)!d_3!d_4!} }$  & 
$\begin{array}{c} d_1+d_3,d_1+d_4\geq 0\\
d_2+d_3, d_2+d_4, d_3,d_4\geq 0\\
(d_1,d_2,d_3,d_4)\neq (0,0,0,0)\end{array}$    \\ \cline{2-4}
& $2,\cdots, 7$ & 0 & \\ \hline

\multirow{7}{*}{$Y_m$} 
& $1,2,3$ & $0$ & \\ \cline{2-4} 
& $4$ & $\begin{array}{c}  \displaystyle{ \frac{(-1)^{2d_1+2d_2-d_3-1}(2d_1+2d_2-d_3-1)!}{(d_1!)^2d_2!(d_{m-1}-2d_m)!d_m!} }\\
\cdot\displaystyle{ \frac{1}{\prod_{b=2}^{m-2} (d_{b-1}-2d_{b}+d_{b+1})!} } \end{array}$ & \\ \cline{2-3}  
& $5,\dots, r-2$ & $\begin{array}{c} \displaystyle{ \frac{(-1)^{d_{i-3}-2d_{i-2}+d_{i-1}-1}}{(d_1!)^2d_2!(-2d_1-2d_2+d_3)!(d_{m-1}-2d_m)!d_m!} }\\
\cdot\displaystyle{ \frac{(-d_{i-3}+2d_{i-2}-d_{i-1}-1)!}{\prod_{b=2,b\neq i-2}^{m-2} (d_{b-1}-2d_{b}+d_{b+1})!} } \end{array}$ & each factorial is non-negative \\ \cline{2-3} 
& $r-1$ & $\begin{array}{c}  \displaystyle{ \frac{(-1)^{2d_m-d_{m-1}-1}(2d_m-d_{m-1}-1)!}{(d_1!)^2d_2!(-2d_1-2d_2+d_3)!d_m!} }\\
\cdot\displaystyle{ \frac{1}{\prod_{b=2}^{m-2} (d_{b-1}-2d_{b}+d_{b+1})!} } \end{array}$ & $(d_1,\dots,d_m)\neq 0$ \\ \cline{2-3}
& $r$ & $\begin{array}{c}  \displaystyle{ \frac{(-1)^{d_m-1}(-d_m-1)!}{(d_1!)^2d_2!(-2d_1-2d_2+d_3)!(d_{m-1}-2d_m)!} }\\
\cdot\displaystyle{\frac{1}{\prod_{b=2}^{m-2} (d_{b-1}-2d_{b}+d_{b+1})!} }  \end{array}$ &  \\ \hline

\end{tabular}
\caption{formulae for $A(q)$ (mirror maps)}
\label{table:A}
\end{table} }
\end{center}

\end{document}